\documentclass[12pt]{amsart}
\usepackage{amssymb,amsfonts,amsmath,amsthm,amscd,array,hhline}
\newtheorem{thm}{Theorem}[section]
\newtheorem{prop}[thm]{Proposition}
\newtheorem{lem}[thm]{Lemma}
\newtheorem{Cor}[thm]{Corollary}

\def\bee{\begin{eqnarray}}
\def\bes{\begin{eqnarray*}}
\def\eee{\end{eqnarray}}
\def\ees{\end{eqnarray*}}
\def\a{\alpha}

\def\s{\sigma}
\def\f{\varphi}
\def\Z{\mathbf Z}
\def\PP{\Phi}
\def\Proof{{\it Proof.\ }\ }
\def\ctd{\hfill$\Box$}
\def\FF{\mathcal F}
\def\VV{\mathcal V}

\begin{document}

\title[Prime (-1,1) and Jordan monsters]
{Prime (-1,1) and Jordan monsters and superalgebras of vector type}											

\author{Sergey V. Pchelintsev}
\address{Finance University under the Government of the Russian Federation,
and
 Moscow City Pedagogical University}
          \email{pchelinzev@mail.ru}
\author{Ivan P. Shestakov}
\address{Instituto de Matem\'atica e Estat\'\i stica,
         Universidade de S\~ao Paulo,
                  S\~ao Paulo, Brazil,
                  and
                  Sobolev Institute of Mathematics, Novosibirsk, Russia}
         \email{shestak@ime.usp.br}

\begin{abstract}

It is proved that the  prime degenerate (-1,1) algebra constructed
in \cite{Pch2}  (the (-1,1)-monster) generates the same variety of
algebras as the Grassman  (-1,1)-algebra. Moreover, the
same variety is generated by the Grassmann envelope of any  simple
nonassociative (-1,1)-superalgebra. The variety occurs to be the
smallest variety of (-1,1)-algebras that contains prime
nonassociative algebras.

Similar results are obtained for Jordan algebras. Thus, the Jordan
monster (the prime degenerate  algebra constructed in \cite{Pch2})
and the Grassmann envelope of the prime Jordan superalgebra of
vector type  have the same ideals of identities. It is also shown
that the Jordan monster generates a minimal variety that contains
prime degenerate Jordan algebras.

All the algebras and superalgebras are considered over a field of
characteristic 0.
\end{abstract}

\maketitle


\section*{Introduction}%

An algebra is called {\em degenerate} if it contains nonzero absolute zero divisors. In \cite{Pch1,Pch2}, the first examples of prime degenerate algebras in the varieties of  alternative, Jordan and algebras of type (-1,1) were constructed, later received the name of {\em Pchelintsev Monsters} \cite{Mc1}.

In \cite{MZ}, there were given another constructions of prime degenerate
Jordan algebras, based on the Jordan superalgebra $J(\Phi [x] ,\frac{d}{dx})$ of vector fields on the line, and the superalgebra of Poisson brackets
$J(\Phi[X,Y];\,[\cdot ,\cdot])$. Observe that the  superalgebras $J(\Phi[x] ,\frac{d}{dx})$ and  $J(\Phi[X,Y];\,[\cdot ,\cdot])$  have different identities.

In \cite{Sh1}, the concept of a (-1,1)-superalgebra of vector type was introduced and it was proved that a Jordan superalgebra of vector type may be obtained as the (super)symmetrized algebra $A^{(+)}$ of a (-1,1)-superalgebra of vector type $A$. Furthermore, the (-1,1)-superalgebras  $A_{VF}:=A(V,\Phi,\tau,\lambda)$ and the Jordan superalgebras $J_{VF}:=J(V,\Phi,\tau,\lambda)$ of vector fields on a line associated with an additive homomorphism $\tau :V\rightarrow \Phi $ of abelian groups and partial map $\lambda :V\rightarrow V$ were introduced in \cite{Sh1}. It was  proved that if $V\neq 0$ and the map $\tau $\ is injective, then the superalgebras $A_{VF}$ and $ J_{VF}$ are
 prime.   These superalgebras were then used to construct prime degenerate (-1,1) and Jordan algebras. 

\smallskip
In all the papers mentioned above, the degenerate prime algebras were constructed as  free algebras in  varieties generated by some auxiliary algebras. Thus, in \cite{Pch2} the auxiliary algebra was
determined by a set of generators and relations, in \cite{MZ} and \cite{Sh1} the auxiliary algebras appeared as the Grassmann envelopes of the superalgebras
$J(\Phi [x] ,\frac{d}{dx})$, $J(\Phi[X,Y];\,[\cdot ,\cdot])$, and $A_{VF}$, $J_{VF},$ respectively.

\smallskip

We give now  the main results of the article, noting that the definitions of all algebras in the statements are given in Sections 2 and~7. Below $G(X)$ denotes the Grassmann envelope of a superalgebra $X$. 

\smallskip

\textbf{Theorem A.} \textit{Let }$A_{0}$\textit{\ be the auxiliary algebra
associated with the prime }(-1,1)\textit{-monster }\cite{Pch2};\textit{\ }$G_{\left(
-1,1\right) }$\textit{\ be the  }(-1,1)\textit{-Grassmann algebra}; $%
A_{VF}$\textit{\ be the }(-1,1)\textit{-superalgebra of vector fields on a
line};\textit{\ }$B$\textit{\ be a simple nonassociative }(-1,1)\textit{%
-superalgebra. Then the algebras }$%
A_{0},\;G_{(-1,1)}$, $G(A_{VF}),\ G(B)$\textit{\ over a field of characteristic }$0$\textit{\ have the same
ideals of identities.}

\smallskip

In particular, we show that the prime (-1,1)-algebras over a field of
characteristic $0$ constructed in \cite{Pch2} and \cite{Sh1} are isomorphic.

\smallskip

We prove Theorem A in  the superalgebra setting, describing the free (-1,1)-superalgebra $\FF_{(-1,1)}[\emptyset;x]$
on one odd generator and its central extensions.
 The case of simple superalgebras is based on the classification of simple nonassociative (-1,1)-superalgebras \cite{Sh2} and on the embedding theorem of a simple nonassociative (-1,1)-superalgebra into a suitable twisted superalgebra of vector type \cite{ZhSh}.

\smallskip

Additional properties of the variety $Var(A_{0}) $ are listed in the following theorem.

\smallskip

\textbf{Theorem B.} \textit{Let }$\mathcal{V}_0=Var\,A_{0}$\textit{\
and let }$\mathcal F(A_{0}) $\textit{\ be the free algebra of
countable rank in the variety }$\mathcal{V}_0$\textit{\ over a field }$\Phi $\textit{\ of characteristic }$0$\textit{. Then}

a) $\mathcal{V}_0$\textit{\ is the smallest variety  of} (-1,1)\textit{-algebras
that contains prime nonassociative algebras};

b)  \textit{The meta-ideals of finite index of the} (-1,1){\it-monster $\mathcal F
(A_{0}) $\textit{\ }have the same identities};

c) \textit{A meta-ideal of the algebra } $\mathcal F(A_{0}) $ \textit{ which is contained in its commutant is not
a free algebra of any variety of algebras.}

\smallskip

Note that in the proof of Theorem A the "minimal" quotient algebra $%
\overline{A}_{0}$, which has the same identities as the algebra $A_{0}$, is constructed. 

It is an open question what identities define the variety $\mathcal{V}_0$? The conjecture  first formulated in \cite{Pch} that $\mathcal{V}_0$ coincides with the variety $\mathcal{S}t$ of {\em strongly (-1,1)} algebras, that is, defined by the identity $[[x,y],z]=0$, is still open.    

\smallskip

The technique developed in the proof of Theorem A, allows us to obtain
certain  analogues of these results for  Jordan algebras.

\smallskip

\textbf{Theorem C.} {\it Let $J=J(\Gamma,\delta)$ be a prime Jordan superalgebra of vector type where $\Gamma=\Gamma_0$, $J[Z;x]$ be a central extension of the free Jordan superalgebra generated by one odd element (see Section 7). Then the algebras}
$$
J_{0}=A_{0}^{(+)},\  G(J),\ (G_{(-1,1)})^{(+)},\ G(J[Z;x]),\ G(J_{VT})
$$
{\it over a field of characteristic 0 have the same ideals of identities.}

\textit{In addition, the variety  }$\mathcal{JV}_0$\textit{\ generated by any of these algebras is a minimal  variety containing prime degenerate algebras}.

\smallskip

As a corollary, it turns out that over a field of characteristic $0$ the prime degenerate special Jordan algebras constructed in \cite{MZ,Sh1} are isomorphic to the Jordan monster $\mathcal F(J_{0}) $. In addition, $\mathcal F(J_{0}) $ coincides with the subalgebra of the Jordan algebra generated by a countable set of free generators $X$ of the algebra $\mathcal F(A_{0}) $.

\smallskip

The simple alternative superalgebras were classified in \cite{ZS, Sh3}. Over a field of characteristic $3$, these are the superalgebras $B(1,2)$, $B(4,2)$, and the twisted superalgebra of vector type $B\left( \Gamma ,D,\gamma \right) $.
The free algebras in the varieties generated by the Grassmann envelopes of finite-dimensional superalgebras $B(1,2)$ and $B(4,2)$ can not be prime (see \cite[page 9]{MZ}),
while for the infinite-dimensional superalgebra $B\left( \Gamma ,D,\gamma \right) $  this question is open.

\smallskip

Concerning the Grassmann algebras, we recall that they may be defined as the algebras of skew-symmetric functions or as certain subalgebras of Grassmann envelopes of free superalgebras generated by one odd element.

For the first time such an algebra was considered by G.V.Dorofeev in the variety of solvable of index two alternative algebras; it is the classical {\em Dorofeev example} \cite{D1}. He used also the nilpotent of index 7 alternative Grassmann algebra to prove that any alternative algebra with three generators satisfies certain identities that do not hold in all alternative algebras \cite{D2}.

\smallskip

In \cite{Sh4,ShZh2},  additive bases for free superalgebras with one odd generator in the varieties of Malcev and alternative superalgebras were constructed, and new elements in the radicals  of free algebras were found. In particular, bases for Malcev and alternative Grassmann algebras were obtained there.

\smallskip

Identities of Grassmann algebras in some varieties of alternative algebras over a field of characteristic $3$ and of right alternative metabelian algebras over a field of characteristic different from $2$ and $3$, have been studied in \cite{Pch5, Pch6, Pch7}.

\smallskip

In this paper we also construct a base of the (-1,1) Grassmann algebra and prove that the free unital (-1,1)-superalgebra with one
odd generator is isomorphic to the (-1,1)-superalgebra  of vector type $B(\Phi[t_0,t_1,\ldots],D,t_0)$ where $D(t_i)=t_{i+1},\,i=0,1,\ldots$.

\smallskip
The basic concepts related with identities and varieties, can be found in \cite{ZSSS}.

\section{The Basic Notions}%
\label{Sec:TheBasicNotions}

\subsection{$(-1,1)-$algebras and Jordan algebras}%
\label{SubSec:algebrasAndJordanAlgebras}

Recall  that an algebra $A$ is called {\em right alternative} if it satisfies the identity 
\bes
(x,y,y)=0,
\ees
where $(x,y,z)=(xy)z-x(yz)$ is the associator of the elements $x,y,z$.
A right alternative algebra $A$ is called a {\em  (-1,1)-algebra} if it
satisfies the identity
\bee
(x,y,z)+(y,z,x)+(z,x,y)&=&0\label{6},
\eee
and it is called a {\em strongly  (-1,1)-algebra} if it satisfies 
\bee
\ [[x,y],z]=0,\label{strong}
\eee
where $[x,y]=xy-yx$ is the commutator of the elements $x,y$.

In any right alternative algebra hold the identities:
\bee
(ab,x,y)+(a,b,[x,y])&=&a(b,x,y)+(a,x,y)b,\label{3}\\
 (a,x,y)x&=&(a,x,xy) \label{4}
\eee
Moreover, any algebra satisfies the identities
\bee
(xy,z,t)-(x,yz,t)+(x,y,zt)=x(y,z,t)+(x,y,z)t,\label{4'}\\ \
[xy,z]-x[y,z]-[x,z]y=(x,y,z)-(x,z,y)+(z,x,y)\label{4''}.
\eee

A commutative algebra is called  {\em Jordan algebra} if it satisfies the identity 
\bes
(x^2,y,x)=0.
\ees

For an algebra $A$ denote by $A^{(+)}$ the associated symmetrized algebra, with  the same vector space and the  multiplication $x\odot y=\frac12(xy+yx)$. If $A$ is a right alternative algebra, then the algebra $A^{(+)}$  is a special Jordan algebra (see \cite{ZSSS}).

By ${A}^{\sharp}$ we will denote the {\em unital hull} of an algebra $A$, that is, the algebra obtained from $A$ by adjoining the external unit
element.

 For  an element $a\in A$ we  denote by $R_{a}$ and $L_a$ the operators of right and left multiplication on $a$:  $R_{a}:x\mapsto xa, \ L_a:x\mapsto ax$.  We use the common notation $T_a$ for any of $R_a, L_a$. We set also $R_{a,b}=R_{a}R_{b}-R_{ab}$.

\subsection{ Superalgebras} Recall that a \textit{superalgebra} $A=A_{0}\oplus A_{1}$ is a  $\mathbb{Z}_{2}$-graded algebra. A typical example of a superalgebra is  the associative Grassmann algebra  $G=G_{0}\oplus G_{1}$, with the generators $1, e_1,e_2,\ldots$;  $e_ie_j=-e_je_i$, and the standard $Z_2$-grading. For a superalgebra
$A=A_{0}\oplus A_{1}$, the {\it Grassmann envelope} $G(A)$ is defined via
$G(A)=A_{0}\otimes G_{0}\oplus A_{1}\otimes G_{1}$.
Let $\mathcal{V}$ be a variety of algebras.  A superalgebra $B=B_{0}\oplus B_{1}$ is called a {\it $\mathcal{V}$-superalgebra} if its Grassmann envelope $G(B)$  is a $\mathcal{V}$-algebra \cite{ZS}. In particular, if $A\in \mathcal{V}$, then $A\otimes G=A_0\otimes G_0\oplus A_1\otimes G_1$ is a $\mathcal{V}$-superalgebra.
We will denote by $F_\mathcal{V}[X;Y]$ the free $\mathcal{V}$-superalgebra on the sets $X$ and $Y$ of even and odd generators, respectively.

It follows easily from the definition that a superalgebra $A=A_{0}\oplus A_{1}$ is a {\em strongly (-1,1)-superalgebra} if and only if for any homogeneous elements $x,y,z$, the following identities hold:
\bes
(x,y,z)+\left( -1\right) ^{yz}(x,z,y)&=&0,\\
\ [[x,y]_{s},z]_{s}&=&0,
\ees
where $\left( -1\right) ^{xy}=\left( -1\right) ^{\left\vert x\right\vert
\left\vert y\right\vert }$, $\left\vert x\right\vert $ means the parity of a homogeneous element $x$, i.e., $\left\vert x\right\vert =i$ if $x\in A_{i}$;
$[x,y]_{s}=xy-\left( -1\right) ^{xy}yx$ is a {\it supercommutator} of  homogeneous elements $x,y$.

\subsection{Centers.} Following A.Thedy \cite{Th}, consider in a right alternative algebra $A$ the following centers:

$K\left( A\right) =\left\{ k\in A\mid \left( \forall x\right) \left[ k,x%
\right] =0\right\} $, {\it the commutative center};

$V\left( A\right) =\left\{ v\in A\mid \left( \forall x\right) \left(
x,x,v\right) =0\right\} $, {\it the left alternative center};

$N\left( A\right) =\left\{ n\in A\mid \left( \forall x,y\right) \left(
x,y,n\right) =\left( x,n,y\right) =\left( n,x,y\right) =0\right\} $, {\it the associative center};

$Z\left( A\right) =K\left( A\right) \cap N\left( A\right) $, {\it the full
center}.

It is known and is easy to see that all the mentioned centers, except $V(A)$, are subalgebras of $A$ and $K\left( A\right) \subseteq V\left(A\right) $.

An important role in the theory of (-1,1)-algebras is played by the commutative center $K(A)$, whose elements we will call \textit{central}. Remind the basic properties of $K(A)$ (see \cite{Hen1}):
\begin{itemize}
\item[a)] $A$ is an associative bimodule over $K(A)$;
\item[b)]  $K(A)$ is invariant under the operators $R_{a,b}$;
\item[c)] if $k\in K(A),\ a,x,y,z\in A$, then
\bee
(k,x,y)&=&2(x,k,y),\label{8}\\
\lbrack kx,y]&=&k[x,y]+\tfrac{3}{2}(k,x,y),\label{9}\\
(x,y,zk)&=&(x,y,z)k+(x,y,k)z,\label{10}\\
(k,a^{2},x)&=&2(k,a,x)a, \label{11}\\
(k^{2},x,y)&=&2(k,x,y)k.\label{12}
\eee
\end{itemize}
Note that identities (\ref{8}) - (\ref{12}) verified also in a (-1,1)-superalgebra $A$ for every $k\in K(A)\cap A_0,\ a\in A_0$, and arbitrary homogeneous elements $x,y,z$.

\subsection{The  Grassmann algebra of a variety $\mathcal{V}$.} Consider
the free superalgebra $F_{\mathcal{V}}[\emptyset ;x]=F_{0}\oplus F_{1}$ in a variety $\mathcal{V}$ on one odd generator $x$.  Following \cite{ShZh1, ShZh2}, we define the $\mathcal{V}$-\textit{Grassmann algebra} $G_{\mathcal{V}}[X]$ as the  subalgebra  of the Grassmann envelope $G\left( F_{\mathcal{V}}[\emptyset ;x]\right) $ generated by the set $X=\{x\otimes e_{i},\ i=1,2,\ldots\}$.

\begin{lem}\label{lem.Grassmann} \cite{ShZh1}
Let $B$ be a base of $F_{\mathcal{V}}[\emptyset;x]$.
Then the $\mathcal{V}$-Grassmann algebra $G_{\mathcal{V}}[X]$ has a base formed by the elements $u(x)\otimes e_{i_1}\ldots e_{i_n}$, where $u(x)\in B,\ n=\deg u, \ i_1<i_2<\cdots i_n$.
\end{lem}

\smallskip
Note that the Grassmann algebra can be defined by generators and relations
as an algebra of skew (alternative) functions. Let $F_{\mathcal{V}}[X]$ be the free
algebra in the variety $\mathcal{V}$ on a set $X=\left\{ x_{i}\mid i=1,2,\ldots\right\} $ of free generators. Consider  the ideal $I$, generated by polynomials of the form: $v, w+w^{\left( ij\right) }$, where $v$ and $w$ are monomials such that $v$ has degree $\geq 2$ with
respect to at least one variable, and $w$ is a multilinear monomial in
the variables $x_{1},\ldots,x_{n}\in X,1\leq i<j\leq n;$ $(ij)$ is a
transposition of  $i$ and $j$. Denote the quotient algebra $F_{\mathcal{V}}[X]/I$ as $Skew_{\mathcal{V}}[X]$ and call it {\em the algebra of skew functions of the variety $\mathcal{V}$}; here the generators $x_{n}$ are identified with their images $x_{n}+I$ under the canonical homomorphism.

\begin{prop} \label{prop1.2} \cite{ShZh1} The algebras $Skew_\mathcal{V}[X]$ and $G_{\mathcal{V}}[X]$ are isomorphic.
\end{prop}

\begin{lem}\label{lem1.3} \cite{Sh4, ShZh1} Let $f$ be a multilinear polynomial of degree $n$. If $f(x,\ldots,x)=0$ in the free $\mathcal{V}$-superalgebra $F_{\mathcal{V}}[\emptyset ;x]$,  then  the identity
$$
\sum\limits_{\sigma \in S_{n}}\left( -1\right) ^{\sigma }f\left(
x_{1\sigma },\ldots,x_{n\sigma }\right) =0
$$
holds in $\mathcal{V}$.
\end{lem}

\begin{Cor}
Let $f$ be a multilinear polynomial of degree $n$. If $f(x_1,\ldots,x_n)=0$ in  $G_{\mathcal{V}}[X]$,  then  the identity
$$
\sum\limits_{\sigma \in S_{n}}\left( -1\right) ^{\sigma }f\left(
x_{1\sigma },\ldots,x_{n\sigma }\right) =0
$$
holds in $\mathcal{V}$.
\end{Cor}

\section{The auxiliary superalgebras}

\subsection{The superalgebras of vector type}
Let $\Gamma =\Gamma_0\oplus\Gamma_1$ be an associative (super)commutative superalgebra,
$D$ be a nonzero even derivation of~$\Gamma$,
$\gamma \in \Gamma_0$.
Denote by $\bar\Gamma$
an isomorphic copy of the vector space~$\Gamma$ with respect to the isomorphism
$a \mapsto \bar a$ and set $ B(\Gamma,D,\gamma) = \Gamma \oplus \bar\Gamma$
with~the multiplication
\begin{align*}
& a \times b = ab,
\\
& a \times \bar b = (-1)^{b}\,(\bar a \times b) = \overline{ab},
\\
& \bar a \times \bar b = (-1)^{b}\,(\gamma ab + 2D(a)b + aD(b)),
\end{align*}
where $a, b \in \Gamma_0\cup\Gamma_1$, $ab$ is the product of the elements $a$~and~$b$
in~$\Gamma$; and with the grading
$$
B(\Gamma,D,\gamma)_{\bar 0} = \Gamma_0+\bar\Gamma_1,\quad
B(\Gamma,D,\gamma)_{\bar 1} = \Gamma_1+\bar\Gamma_0.
$$
The superalgebra $B(\Gamma,D,\gamma)$ is called {\it the twisted superalgebra of vector type} \cite{Sh1}. It is a strongly (-1,1)-superalgebra which is simple if and only if the algebra~$\Gamma$ does not contain proper  $D$-invariant
ideals (i.e. is $D$-simple).

 The adjoint symmetrized superalgebra $B(\Gamma,D,\gamma )^{(+)}$ with the supersymmetric multiplication
$$
x\bullet y=\frac{1}{2}\left( x\times y+(-1)^{xy}y\times x\right)
$$
is isomorphic to the Jordan algebra of vector type $J(\Gamma,\delta)$, i.e., has the following multiplication:
$$
a\bullet b=ab,\qquad a\bullet \overline{b}=(-1)^ b\,\overline{a}\bullet b=\overline{%
ab},\qquad \overline{a}\bullet \overline{b}=(-1)^ b(a^{\delta }b-ab^{\delta }),
$$
where $a^{\delta }=\frac{1}{2}D\left( a\right) $. The superalgebras $J(\Gamma ,\delta )$ were introduced by K.McCrimmon \cite{Mc}.

\smallskip

When $\Gamma=\Gamma_0$ is an algebra, we will call the superalgebras $B(\Gamma,D,\gamma)$ and $J(\Gamma,\delta)$  of {\em  even vector type}.

\subsection{The superalgebras $A_{VF}$ and $J_{VF}$ of vector fields on a line} In \cite{Sh1}, the second author introduced  (-1,1)-superalgebras of vector fields on a line $A_{VF}:=A\left( V,\Phi ,\tau ,\lambda \right) $.
Recall that $V$ is an additive commutative semigroup; $\tau :V\rightarrow
\Phi $ is an additive homomorphism; $\lambda :V\rightarrow V$ is a partial
map defined everywhere except possibly the neutral element 0 that satisfies the condition $\lambda \left( u+v\right) =\lambda \left( u\right) +v$ in its domain.

Consider the graded vector space $A=A_{0}\oplus A_{1}$ with the bases
$\left\{ a_{v}\mid v\in V\right\} $ and $\left\{ x_{v}\mid v\in V\right\} $
for $A_{0}$\ and $A_{1}$, respectively. The multiplication on the superalgebra $A_{VF}=A$ is defined  by the rules:
$$
a_{u}\cdot a_{v}=a_{u+v},\  a_{u}\cdot x_{v}=x_{u}\cdot
a_{v}=x_{u+v},\ x_{u}\cdot x_{v}=\left( 4\tau \left( u\right) +2\tau
\left( v\right) \right) a_{\lambda \left( u+v\right) }.
$$

The Jordan superalgebra $J_{VF}$ is obtained from $A_{VF}$ by supersymmetrization: $J_{VF}=A_{VF}^{(+)}$. Each of the superalgebras $A_{VF}$ and $J_{VF}$ is a superalgebra of even vector type. Also, if $V\neq 0$ and $\tau$ is injective, then the superalgebras $A_{VF}$ and $J_{VF}$ are prime \cite[Theorem 1]{Sh1}.

\subsection{The algebra $A_{0}$.} Following \cite{Pch2}, we denote by  $A_{0}$ the strongly (-1,1)-algebra defined by the generators $z,e_{1},\ldots ,e_{n},\ldots $ and the relations:
$$
e_{i}e_{j}+e_{j}e_{i}=0,\quad
(e_{i},e_{j},e_{p})=0,\quad
([e_{i},e_{j}],e_{p},e_{q})=0,
$$
$$
((z,e_{i},e_{j}),e_{p},e_{q})=0, \qquad \ \ \
[z,f(z,e_{i},\ldots ,e_{j})]=0
$$
 for any polynomial $f(z,e_{i},\ldots ,e_{j})$.
It is known \cite{Pch2} that the algebra $A_{0}$ has an additive basis consisting of the elements:
\bee
z^{n}(z,g_{1},g_{2})\cdots (z,g_{2m-1},g_{2m})(h_{1}\cdots h_{p}),
\eee
where $n,m,p$ are nonnegative integers with $n+m+p>0$, $g_{i},h_{i}$ are elements of the set $E=\left\{ e_{1},\ldots ,e_{n},\ldots
\right\} $, which is ordered by its indices and $g_{1}<\cdots
<g_{2m}<h_{1}<\cdots <h_{p}$. Observe that the elements of the form $(z,e,g)$ and $eg$, where $e,g\in E$, are contained in the full center $Z\left( A_{0}\right)$ of the algebra $A_{0}$. Note also that the subalgebra generated by the set $E$ is an associative Grassmann algebra (without the identity) with the standard generators $e_{1},\ldots ,e_{n},\ldots $.

\smallskip
The algebra $A_{0}$ can be constructed through superalgebras. Consider the   twisted superalgebra  of even vector type $B=B\left(\Phi \left[ z,s\right],s\tfrac{d}{dz},1\right)$. Observe that
$$
\bar 1\times \bar 1=1,\ (z,\bar 1,\bar 1)=2s,\ (s,\bar 1,\bar 1)=0.
$$
  Now the subalgebra of the Grassmann envelope $G(B)$ generated by the elements $z\otimes 1,\ \bar 1\otimes e_{i},\ i=1,2,\ldots,$ is isomorphic to the algebra $A_{0}$.

\subsection{The auxiliary algebras $\overline{A}_{0}$ and $\overline{J}_{0}$.} Denote by $\overline{A}_{0}$ the algebra
defined in the variety of strongly (-1,1)-algebras by the generators $z,e_{1},\ldots ,e_{n},\ldots $ and the relations:
$$
e_{i}e_{j}=0, \qquad ((z,e_{i},e_{j}),e_{p},e_{q})=0, \qquad \lbrack
z,f(z,e_{i},\ldots ,e_{j})]=0
$$
for any polynomial $f(z,e_{i},\ldots ,e_{j}) $.
It is known \cite{Pch2} that the algebra $\overline{A}_{0}$ has an additive basis consisting of the elements:
\bee\label{14}
z^{n}(z,g_{1},g_{2})\ldots \ (z,g_{2m-1},g_{2m})g_{2m+1},
\eee
where $n$ and $m$ are nonnegative integers, $n+m>0$, $g_{i}\in E=\left\{ e_{1},\ldots ,e_{n},\ldots \right\}$ and $g_{1}<\cdots \ <g_{2m}<g_{2m+1}$. Note that the elements of the form $(z,e,g) $, where $e,g\in E$ are contained in the full center of the
algebra $Z\left( \overline{A}_{0}\right) $.

\smallskip
Similarly to $A_0$, the algebra $\overline{A}_{0}$ can be constructed via superalgebras, using the twisted superalgebra of vector type $B\left(\Phi \left[z,s\right],s\tfrac{d}{dz},0\right) $.

\smallskip

By analogy with the $\overline{A}_{0}$, denote by $\overline{J}_{0}$ the algebra defined in the variety of Jordan algebras by the generators $z,e_{1},\ldots ,e_{n},\ldots $ and the relations:
$$
e_{i}e_{j}=0,\ ((z,e_{i},e_{j}),e_{p},e_{q})=0,\ (z,f,f)=0
$$
 for any $f=f(z,e_{i},\ldots ,e_{j})$.
It is clear that the algebra $\overline{J}_{0}$ has an additive basis
consisting of the elements of the form (\ref{14}). In addition, it is easy to see that

a) $(z,e,g)\in Z\left( \overline{J}_{0}\right) $ for any $e,g\in E$;

b) $\overline{J}_{0}=\overline{A}_{0}^{(+)}$.

\section{ The free (-1,1)-superalgebra $F_{(-1,1)}
[\emptyset;x]$ and its central extension $F[Z;x]$}

\subsection{The free (-1,1)-superalgebra $F_{(-1,1)}[\emptyset;x]$} Let $\Phi[T]$ be the algebra of
polynomials on a countable set of variables
$T=\{t_0,\ldots,t_n,\ldots\}$, and let $D$ be the derivation of
$\Phi[T]$ defined by the condition $D(t_i)= t_{i+1},\
i=0,1,\ldots.$ Consider the twisted superalgebra of vector type
$B(\Phi[T],D,t_0)$. By \cite[Theorem 4]{Sh3}, the superalgebra
$B(\Phi[T],D,t_0)$ is prime. Let $\Phi_0[T]$ denotes the
subalgebra of polynomials without constant terms, then the
subspace $\Phi_0[T]\oplus\overline{\Phi[T]}$ is a subsuperalgebra
of $B(\Phi[T],D,t_0)$ which we will denote as
$B_0(\Phi[T],D,t_0)$. It is clear that
$B(\Phi[T],D,t_0)=(B_0(\Phi[T],D,t_0))^{\sharp}$, the unital hull.
Clearly, $B_0(\Phi[T],D,t_0)$ is also prime. One can easily check
that it is generated by the odd element $\bar 1$.

\begin{thm}\label{thm3.1}
Let $\Phi$ be a field of characteristic $\neq 2,3$. The free
$(-1,1)$-superalgebra $F_{(-1,1)}[\emptyset;x]$ over $\Phi$ on one
odd generator $x$ is isomorphic to the superalgebra
$B_0(\Phi[T],D,t_0)$ with the free generator $\bar 1$.
\end{thm}

Denote $F:=F_{(-1,1)}[\emptyset;x]$. We prove first two lemmas.

\begin{lem}\label{lem3.2}
 $x^2\in K(F)$.
\end{lem}
\Proof By superized associator Jacoby identity (\ref{6}), we have
$$
3(x,x,x)=0,\  (x^2,x,x)-(x,x^2,x)+(x,x,x^2)=0,
$$
which implies
\bee
(x,x,x)&=&[x^2,x]=0,\label{(xxx)}\\
(x^2,x,x)&=&2(x,x^2,x).\label{2xxx}
\eee
 By  associator identity (\ref{4'}), we have
\bes
(x^2\cdot x,x,x)-(x^2,x^2,x)+(x^2,x,x^2)&=&x^2(x,x,x)+(x^2,x,x)x,\\
(x^2,x,x^2)-(x,x^2,x^2)+(x,x,x\cdot x^2)&=&x(x,x,x^2)+(x,x,x)x^2,\\
(x^2,x^2,x)-(x,x\cdot x^2,x)+(x,x,x^2\cdot x)&=&x(x,x^2,x)+(x,x,x^2)x,
\ees
which implies
\bes
(x^3,x,x)&=&2(x^2,x^2,x)+(x^2,x,x)x,\\
(x,x,x^3)&=&(x^2,x^2,x)+x(x,x,x^2),\\
(x^2,x^2,x)&=&[x,(x,x^2,x)].
\ees
Now, the superized associator Jacoby identity (\ref{6}) implies
$$
(x^3,x,x)+2(x,x,x^3)=0,
$$
hence by (\ref{2xxx}) we have  $4(x^2,x^2,x)+2[(x,x^2,x),x]=0$ and eventually
\bee\label{2x2xx}
(x^2,x^2,x)=0.
\eee
In particular, $[x^3,x^2]=(x^2,x,x^2)=0$.

Let now $u(x)$ be a monomial on $x$ of degree $n>3$, consider $[u(x),x^2]$. We may assume that $[v(x),x^2]=0$ for every monomial $v(x)$ of degree less then $n$.
 Recall that a (-1,1)-algebra $A$ satisfies the identity $[[a,b],D(A)]=0$ where $D(A)$ is the associator ideal of $A$ (see \cite{Roo1}).
Therefore, in $F$ we have $[2x^2,D(F)]=[[x,x]_s,D(F)]=0$, and hence we may write
\bes
[u(x),x^2]=[x^2(x^2v(x)),x^2],
\ees
where $\deg v(x)\geq 0$. Denote $a=x^2,\ v=v(x)$. By (\ref{4''}) and the induction assumption we have
\bes
[u(x),x^2]&=&[a(av),a]=a[av,a]+[a,a](av)\\
&&+(a,av,a)-(a,a,av)+(a,a,av)\\
&=&-(a,a,av). \ees Since $a=x^2$ is even, by (\ref{4}) we have
$(a,a,av)=(a,a,v)a$. The above arguments show that
$(a,a,v)=-[av,a]=0$, proving the lemma. \ctd
\begin{lem}\label{lem3.3}
$[[F,F]_s,F]_s=0$, that is, $F$ is a strongly (-1,1)-super\-al\-ge\-bra.
\end{lem}
\Proof
Evidently, it suffices to prove that $u(x)\in K(F)$ for every monomial $u$ of even degree. We will use induction on $\deg u$. A base of the  induction is given by Lemma \ref{lem3.2}.

Since the variety of (-1,1)-algebras is a $2$-variety (the square of an
ideal is an ideal), we may assume that $u$ has form  $vx$ or $xv$, where $v$ is a monomial of odd degree  (see, for instance,  \cite{ZheShe}).
Similarly, we may assume that $v=wx$ or $v=xw$ for some even $w\in F$.
The inclusion $x^2\in K(F)$ implies, due to \cite[theorem 13.10]{ZSSS}, that
\bee
yz+zy\in K(F) \hbox{ for any odd } y,z\in F.\label{yz+zy}
\eee
In particular, $vx+xv\in K(F)$, and it suffices to prove that $(wx)x,\ x(xw)\in K(F)$. By induction, $w\in K(F)$, hence $wx^2, x^2w\in K(F)$, and  we have by (\ref{9})
\bes
(wx)x&=&wx^2+(w,x,x)\equiv (w,x,x)=\tfrac{2}{3}([wx,x]_s-w[x,x]_s)\\
&=&\tfrac{2}{3}((wx)x+x(wx)-2wx^2)\stackrel{(\ref{yz+zy})}\equiv  0\ (mod\ K(F))
\ees
Similarly, by (\ref{8}),
$$
x(xw)=x^2w-(x,x,w)\equiv -(x,x,w)\equiv \tfrac12(w,x,x)\equiv 0\ (mod\ K(F)).
$$
\ctd

{\it Proof of the theorem.} Denote $A=F_0$, then $F_1=A^{\sharp}x,\
F=A\oplus A^{\sharp}x$, where by  1.4.a) $A$ is an associative and
commutative $\Phi$-algebra and $A^{\sharp}x$ is a commutative and
associative $A$-module generated by $x$. Furthermore, by
(\ref{12}), the application $D:a\mapsto \tfrac12 (a,x,x)$ is a
derivation of $A$. We claim that $A$ coincides with the
$\Phi$-subalgebra $A_0$  generated by the set $\{x^2, D^i(x^2),\
i=1,2,\ldots\}$. It is equivalent to say that $F=A_0+A_0^{\sharp}x$, and since $x\in 
A_0+A_0^{\sharp}x$, it suffices to prove that $A_0+A_0^{\sharp}x$ is a subsuperalgebra of $F$. Since $A_0\subseteq K(F)$, we have $a\cdot bx=bx\cdot a=(ab)x$ for any $a,b\in A$. 

Now, let us calculate the product of elements from $A^{\sharp}x$:
 \bes
ax\cdot bx&=&(a,x,bx)+a(x\cdot bx)=(a,bx,x)+a(x\cdot xb)\\
&=&(ab\cdot x)x-a(bx\cdot x)+a\cdot x^2b-a(x,x,b)\\
&=&(ab,x,x)+(ab)x^2-a(b,x,x)-a(bx^2)+(ab)x^2-a(x,x,b)\\
&=&(\hbox{by } (\ref{8})) = 2D(ab)-2aD(b)+ aD(b)+(ab)x^2\\
&=& 2D(a)b+aD(b)+(ab)x^2.
 \ees
Therefore, $A_0+A_0^{\sharp}x$ is a subsuperalgebra of $F$ and $F=A_0+A_0^{\sharp}x, \ A_0=A=F_0$. Moreover, the obtained relation shows that the superalgebra $F$ is a homomorphic image of the  superalgebra $B_0(A^{\sharp},D,x^2)$ under the homomorphism $\pi: a+\bar b\mapsto a+xb$. Consider the homomorphism $\varphi
:\Phi_0[T]\rightarrow A,\ t_i\mapsto D^i(x^2)$. Clearly, it is a
homomorphism of differential algebras which can be extended to a
homomorphism of superalgebras
$$
\tilde\varphi: B_0(\Phi[T], D,t_0)\rightarrow  B_0(A^{\sharp},D,x^2).
$$
Now the composition $\pi\circ\tilde\varphi$ maps surjectively
$B_0(\Phi[T], D,t_0)$  onto $F$. Since $F$ is a free superalgebra
generated by $x$, this map is an isomorphism. \ctd

\begin{Cor}\label{cor3.4}
 $(F_{(-1,1)}
[\emptyset;x])^{\sharp}\cong B(\Phi[T],D,t_0)$.
\end{Cor}
\begin{Cor}\label{cor3.5}
The variety $Var\left( F_{(-1,1)}\left[ \emptyset ;x\right]
\right)$ can not be generated by a finite dimensional superalgebra.
\end{Cor}
\Proof
 Suppose by contradiction that $Var\left( F_{\left( -1,1\right) }^{s}\left[
\emptyset ;x\right] \right) =Var\left( B\right) $ where $\dim
B=m<\infty $. Consider the Grassmann envelope $G(B)$.
 Let $x_1,\ldots,x_k,\ldots,x_{m}$ be a base of $B$ with $x_1,\ldots,x_k$ be a base of $B_0$, then any element $u\in [G(B),G(B)]$ can be
written in the form $u=\sum_{i,j=1}^ka_{ij}\otimes  [x_i,x_j]+\sum_{i,j=1}^{m} g_ig_j\otimes x_i x_j$, where $a_{ij}\in G_0,\, g_i\in G_1$  or $g_j\in G_1$. 
Clearly,  $(g_ig_j)^2=0$ for every pair $i,j$. Besides, by the identity of Kleinfeld $[x,y]^3=0$ valid in every strongly (-1,1)-algebra \cite{Kl},
we have $[x_i,x_j]^3=0$ for $i,j\leq k$. Since $[G(B),G(B)]\subseteq Z(G(B))$, we see that the element $u$ is nilpotent of degree $m^3+1$.  In other words, the subspace $[G(B),G(B)]\subseteq Z(G(B))$ satisfies the identity $x^{m^3+1}=0$.
Linearizing this identity, in view of associativity and
commutativity of  $Z(G(B))$ we get that $[G(B),G(B)]$
satisfies the identity $(m^3+1)! \,x_1x_2\ldots x_{m^3+1}=0$.
Since $char\,\Phi=0$, this implies that $G(B)$ satisfies the
identity
$$
[x_1,y_1][x_2,y_2]\cdots [x_{m^3+1},y_{m^3+1}]=0.
$$
Therefore, the superalgebra $B$ satisfies the identity
$$
[x_1,y_1]_s[x_2,y_2]_s\cdots [x_{m^3+1},y_{m^3+1}]_s=0.
$$
Since $F_{(-1,1)}\left[ \emptyset ;x\right]\in Var\,B$, it should
satisfies this identity as well, which is impossible since
$[x,x]_s^k=(2x^2)^k\neq 0$ for any $k$. \ctd

\subsection{The central extension $F[Z;x]$}
Let $F_{St}[Z;x]$ be the free strongly (-1,1)-superalgebra on a set $Z$ of even generators and on an odd generator $x$. Denote by $F[Z;x]$ the quotient superalgebra $F_{St}[Z;x]/I$ where $I$ is the ideal generated by the set $\{[z,f],\, z\in Z,\ f\in F_{St}[Z;x]\}$. In other words, $F[Z;x]$ is a strongly (-1,1)-superalgebra freely generated by a set of even central variables $Z$ and an odd variable $x$. Clearly, $F[Z;x]$ satisfies the following universal property:\\
 {\it for any strongly
{(-1,1)}-superalgebra $B$, any odd $b\in B_1$,  and for any
mapping $\varphi :Z\rightarrow K(B)$ there exists a unique
homomorphism $\tilde\varphi:F[Z;x]\rightarrow B$ such that
$\tilde\varphi|_Z=\varphi,\ \tilde\varphi(x)=b.$}

\smallskip

Consider the polynomial ring $\Phi[T_Z]$, where 
$T_Z=T_0\cup(\cup_{z\in Z}T_z)$, $T_0=\{t_0,t_1,\ldots\},
T_z=\{z=z_0,z_1,z_2,\ldots\},\ z\in Z$, and let $D$ be the
derivation of $\Phi[T_Z]$ defined by $D(t_i)=t_{i+1},\
D(z_i)=z_{i+1}$.
\begin{prop}\label{prop3.6} Let $\Phi$ be a field of characteristic $\neq 2,3$. The
superalgebra $F[Z;x]$ is isomorphic to the superalgebra $B_0(\Phi[T_Z], D,t_0)$.
\end{prop}
\Proof  Denote $F=F[Z;x],\  K=K(F)$, and let $A$ be a $\Phi$-subalgebra of $F_0$ generated by the set $\{D^i(x^2),D^j(z)\,|\,z\in Z; i,j=0,1,\ldots\}$, where $D(a)=\tfrac12 (a,x,x)$. Let us prove that $F=A+A^{\sharp}x$.
Observe first that since $F$ is strongly (-1,1), we have $0=[F,[x,x]_s]=2[F,x^2]$, hence $x^2\in K$.  Furthermore, by (\ref{12}) $D$ is a derivation of $K$, hence 
 $A\subseteq K$. Since $Z\cup\{x\}\subseteq A+A^{\sharp}x$, it suffices to prove   
that $A+A^{\sharp}x$ is a subsuperalgebra of $F$.

As above, we have for any $a,b\in A^{\sharp}$
$$
a\cdot bx=ab\cdot x,\ ax\cdot bx= 2D(a)b+ aD(b)+(ab)x^2.
$$
This proves that $F=A+A^{\sharp}x$.  Moreover, this also shows that 
 the superalgebra $F$ is a homomorphic image of the
superalgebra $B_0(A^{\sharp}, D,x^2)$ under the homomorphism $\pi: a+\bar b\mapsto a+xb$. Consider the homomorphism $\varphi
:\Phi_0[T_Z]\rightarrow A,\ t_i\mapsto D^i(x^2),\ z_i\mapsto
D^i(z),\,z\in Z$. Clearly, it is a homomorphism of differential
algebras which can be extended to a homomorphism of superalgebras
$$
\tilde\varphi: B_0(\Phi[T_Z], D,t_0)\rightarrow B_0(A^{\sharp}, D,x^2).
$$
Now the composition $\pi\circ\tilde\varphi$ maps surjectively
$B_0(\Phi[T_Z], D,t_0)$  onto $F$. By the universal property of $F$, this map is invertible and hence it is an isomorphism.
 \ctd
\begin{Cor}\label{cor3.7}
The superalgebra $F[Z;x]$ has a  base of the form $B\cup Bx$,
where $B$ is the set of commutative and associative monomials on
the variables
$$
x^2R_{x,x}^i,\ zR_{x,x}^j,\ z\in Z; \ i,j=0,1,\ldots.
$$
\end{Cor}

The proof follows easily from the isomorphism given in the
theorem.

\ctd

\smallskip
Consider the superalgebra $F[z;x]\cong B_0(\PP[T_0\cup
T_z],D,t_0)$. We will be interested in the following homomorphic
images of $F[z;x]$:
\begin{itemize}
\item  $F_{t}=F[z;x]/((x^2,x,x), ((z,x,x),x,x)\cong
B_0(\PP[t,z,s],s\tfrac{d}{dz},t),$
\item  $F_1=F[z;x]/(x^2-1, ((z,x,x),x,x))\cong
B(\PP[z,s],s\tfrac{d}{dz},1),$
\item  $F_0=F[z;x]/(x^2, ((z,x,x),x,x))\cong
B_0(\PP[z,s],s\tfrac{d}{dz},0)$.
\end{itemize}

Let $A=A_1\oplus A_1$ be a (-1,1)-superalgebra, $a\in A_0,\ x\in A_1$.
Let, furthermore,  $\Z_+$ be the set of non-negative integers and $\Z_+^{\infty}=\cup_n(\Z_+)^n$ be the set of all ordered finite sequences of elements of $\Z_+$.  Denote, for  $I=(i_0,i_1,\ldots,i_k)\in \Z_+^{\infty}$,  
$$
a^I=a^{i_0}(aR_{x,x})^{i_1}(aR_{x,x}^2)^{i_2}\cdots (aR_{x,x}^k)^{i_k}.
$$
Set also $|I|=i_0+\cdots +i_k$ and $d(I)\leq k+1$ to be the number of nonzero elements in $I$. 
 
\begin{Cor}\label{cor3.8}
The superalgebras $F[z;x],\ F_t,\ F_1,\ F_0$ have following bases
in terms of generators $z,x$:
 \bes
F[z;x]&:& (x^2)^Iz^Jx^{\varepsilon},\ I, J\in \Z_+^{\infty},\,  
 \varepsilon\in\{0,1\},\ \varepsilon+|I|+|J|>0;\\
 F_{t}&:&  (x^2)^m z^{n}(zR_{x,x})^{k}x^{\varepsilon},\  m,n,k\geq 0,\, \varepsilon\in\{0,1\},\, \varepsilon+m+n+k>0; \\
 F_1 &:& \ z^{n}(zR_{x,x})^{k}x^{\varepsilon},\ n,k\geq 0,\,
 \varepsilon\in\{0,1\};\\
F_0 &:& z^{n}(zR_{x,x})^{k}x^{\varepsilon},\ n,k\geq 0,\,
\varepsilon\in\{0,1\},\, \varepsilon+n+k>0.
 \ees
\end{Cor}
\begin{Cor}\label{cor3.9}
The Grassmann (-1,1)-algebra $G_{(-1,1)}$ has a following base 
$$
(x^2)^Ix^{\varepsilon}\otimes e_{j_1}e_{j_2}\cdots e_{j_n},\ j_1<j_2<\cdots<j_n,\ I\in \Z_+^{\infty},\ n=\deg ((x^2)^Ix^{\varepsilon}).
$$
\end{Cor}
The proof follows from the previous Corollary and Lemma \ref{lem.Grassmann}. 
\ctd

\section{Small varieties of superalgebras}

We call a variety $\mathcal{V}$ of strongly (-1,1)-superalgebras to be {\it
small} if it does not contain the superalgebra $F_0$.

The following result gives a criterion for a unitary closed
variety $\mathcal{V}$ to be small.

\begin{prop}\label{prop4.1}
A  variety $\mathcal{V}$ of strongly
(-1,1)-su\-per\-al\-gebras over a field $\PP$ of characteristic 0
is small if and only if it satisfies the identity
 \bee\label{id4.1}
(z,x,x)^k=0
 \eee
for any central $z$ and odd $x$ and some $k>0$.
\end{prop}

 \Proof Note that for any
$n$ the element $(z,x,x)^n$ belongs to the base of $F_0$ given in
Corollary \ref{cor3.8}, hence it is non-zero. Therefore, the
condition above is sufficient for a variety $\mathcal{V}$ to be small.

Assume now that $F_0\not\in \mathcal{V}$. Consider the free $\mathcal{V}$-superalgebra
$F_\mathcal{V}[a;x]$ on an even generator $a$ and odd generator $x$. Then
the quotient superalgebra $S=F_\mathcal{V}[a;x]/(x^2, ((a,x,x),x,x),
[a,F_\mathcal{V}[a;x]])$ is a proper homomorphic image of $F_0$. Therefore,
the images of basic elements of $F_0$ are linearly dependent in
$S$. In other words, we have in $F_\mathcal{V}[a;x]$ a relation of the form
 \bes
\sum_{n,k}\a_{n,k}a^{n}(a,x,x)^{k}x^{\varepsilon}=x^2
U+((a,x,x),x,x)V+\sum_i[a,f_i]W_i
 \ees
for some $f_i=f_i(a,x)\in F_\mathcal{V}[a;x]$ and
$U=U(a,x),V=V(a.x),W_i=W_i(a,x)$ from the multiplication algebra
of $F_\mathcal{V}[a;x]$. Since $\PP$ is infinite, we may assume that the
relation above is homogeneous in $a,x$, that is, has a form
 \bes
a^{n}(a,x,x)^{k}x^{\varepsilon}=x^2
U+((a,x,x),x,x)V+\sum_i[a,f_i]W_i.
 \ees
The application $D=R_{x,x}$ in view of (\ref{3}) is a derivation  of the superalgebra $S$ such that $D(x)=D((a,x,x))=0$, hence we have in $S$ 
\bes
D^n(a^{n}(a,x,x)^{k}x^{\varepsilon}) = (a,x,x)^{k+n}x^{\varepsilon}.
\ees
Returning to  $F_\mathcal{V}[a;x]$,   we will get a relation of the form
 \bes
(a,x,x)^{k}x^{\varepsilon}=x^2 U+((a,x,x),x,x)V+\sum_i[a,f_i]W_i.
 \ees
If $\varepsilon=1$, multiplying the relation above by $a$ and $x$
we get in the superalgebra $S$
 \bes
  (((a,x,x)^{k}x)a)x&=&((a,x,x)^{k}a)x\cdot
  x=((a,x,x)^{k}a,x,x)\\
  &=&((a,x,x)^{k},x,x)a+(a,x,x)^{k}(a,x,x)\\
  &=&k(a,x,x)^{k-1}((a,x,x),x,x)a+(a,x,x)^{k+1}\\
  &=& (a,x,x)^{k+1}.
 \ees
Therefore, without loss of generality, we may assume that we have
a homogeneous relation in $F_\mathcal{V}[a;x]$ of the form
 \bes
  (a,x,x)^{k}=x^2U+((a,x,x),x,x)V+\sum_i[a,f_i]W_i.
 \ees
Denote $F_\mathcal{V}[z;x]=F_\mathcal{V}[a;x]/([a,F_\mathcal{V}[a;x]])$, then we have in $F_\mathcal{V}[z;x]$:
 \bes
(z,x,x)^{k}=x^2U(z,x)+((z,x,x),x,x)V(z,x).
 \ees

Consider in the free central extension  $F[s,z;x]$ the element
$u=sU(z,x)$. By Corollary \ref{cor3.7}, $u$ may be written as a
linear combination of monomials on the variables
 $x^2R_{x,x}^i,\,sR_{x,x}^j,\,zR_{x,x}^l$. By homogeneity, since
 $\deg_z(u)=k$, every monomial contains $k$ variables of type $zR_{x,x}^l$.
Since $\deg_xu=2(k-1)$, at least one of these variables should be
just $z$, henceforth $u=zu_1$.  Since $x^2$ is a central element
in $F_\mathcal{V}[z;x]$, we have a homomorphism $\f : F[s,z;x]\rightarrow F_\mathcal{V}[z;x],\ x\mapsto x,\, z\mapsto z,\, s\mapsto x^2$.  Then $x^2U(z,x)=\f(u)=z\f(u_1)$.


Similarly, the element $v=sV(z,x)\in F[s,z;x]$ may be written as
$v=zv_1$ and therefore $((z,x,x),x,x)V(z,x)=z \f(v_1)$ where $\f:
F[s,z;x]\rightarrow F_\mathcal{V}[z;x],\ x\mapsto x,\, z\mapsto z,\,
s\mapsto ((z,x,x)x,x)$.

\smallskip
Resuming, we have in $F_\mathcal{V}[z;x]$ a relation
 \bes
(z,x,x)^{k}=zw(z,x).
 \ees
Linearizing this relation on $z$, we get  
$$
(z_1,x,x)\cdots (z_k,x,x)=\sum_{i=1}^k z_iw_i.
$$ 
Substitute $z_1^2$ instead of $z_1$:
$$
(z_1^2,x,x)\prod_{i\geq 2} (z_i,x,x)=z_1^2w_1+\sum_{i\geq 2} z_iw'_i.
$$ 
But $(z_1^2,x,x)=2z_1(z_1,x,x)$, hence
$$
(z_1^2,x,x)\prod_{i\geq 2} (z_i,x,x)=2z_1(z_1,x,x)\prod_{i\geq 2} (z_i,x,x)=2z_1^2w_1+\sum_{i\geq 2}z_iw''_i.
$$
Comparing the two expressions, we get
$
z_1^2w_1=\sum_{i\geq 2}z_i(w'_i-w_i''),
$
which implies 
$$
(z_1^2,x,x)\prod_{i\geq 2} (z_i,x,x)=\sum_{i\geq 2} z_iw^{(1)}_i.
$$ 
Repeating the same arguments, we get after $k-1$ steps
$$
\left(\prod_{i=1}^{k-1}(z^2_i,x,x)\right)(z_k,x,x)= z_kw_k^{(k-1)},
$$ 
and eventually 
$$
\prod_{i=1}^{k}(z^2_i,x,x)=0.
$$
Due to the inclusion $D(A)^3\subseteq D(A^2)A^{\sharp}$ which holds for any algebra $A$ with a derivation $D$, we have for $D=R_{x,x}$
$$
(z,x,x)^3=\sum_i (z_i^2,x,x)z_i'\ \hbox{ for some } z_i,z_i'\in K(F_\mathcal{V}[z;x]),
$$
and finally $(z,x,x)^{3k}=0$.
\ctd 

\section{The superalgebras $F_0,\, F_1,\, F_{(-1,1)}[\emptyset;x]$ have  the same identities}
We will start with the following auxiliary result.

\begin{lem}\label{lem5.0}
Let $A$ be a superalgebra, $Z=Z(A)$ be the (full) cener of $A$. Assume that for any $a\in A$ there exists $z\in Z$ such that $az\neq 0$. Then the variety $Var\,A$ is unitary closed.
\end{lem}
\Proof
Recall that a variety $V$ is unitary closed if for any element $f$ from the $T$-ideal $T(V)$ of identities of $V$ all the partial derivatives $f\Delta_i$ belong to $T(V)$  as well (see \cite{ZSSS}). Since $char\,\Phi=0$, we may assume that $f$ is multilinear. Assume that $Var\,A$ is not unitary closed, then there exists a multilinear $f=f(x_1,\ldots,x_n)\in T(A)$ such that $f\Delta_i\notin T(A)$ for some $i$. Let, for example,  $f\Delta_n\notin T(A)$, then there exist  $a_1,\ldots,a_{n-1}\in A$ such that  $(f\Delta_n)(a_1,\ldots,a_{n-1})=b\neq 0$. Choose $z\in Z$ with $bz\neq 0$, then we have 
$$
f(a_1,\ldots,a_{n-1},z)= (f\Delta_n)(a_1,\ldots,a_{n-1})z=bz\neq 0.
$$
The contradiction proves the Lemma. \ctd

\begin{lem}\label{lem5.1}
Let $F_0[Z;x]=F[Z;x]/(x^2)$. Then the superalgebras $F_0[Z;x]$ and
$F_0$ generate the same variety.
\end{lem}
\Proof Since $F_0$ is a homomorphic image of $F_0[Z;x]$, it
suffices to prove that $F_0[Z;x]$ belongs to the variety $\mathcal{V}$
generated by $F_0$. Consider first the case when $Z$ is a
singleton: $Z=\{z\}$. Assume that $F_0[z;x]\not\in \mathcal{V}$, then some
nontrivial linear combination $f=f(z,x)$ of basic elements of
$F_0[z;x]$ vanishes in $F_0$ for any substitution $z=a\in
K(F_0),\,x=y\in(F_0)_1$. More exactly, we may assume that
 $f=\sum_{I}\a_{I}z^Ix$, where $I=(i_0,\ldots,i_k)\in \Z_+^{\infty},\, \a_{I}\in\PP$.
  Since $char\,\PP=0$, we may assume that $f$ is
homogeneous, that is, the numbers $|I|$ and $i_1+\cdots
+ki_k=i$ are fixed. Let $m$ be a maximal value of the index $i_0$ for  $I$ in $f$. The algebra $F_0$ contains in its center the polynomial ring $\PP[s]$ which acts without torsion on $F_0$. By Lemma \ref{lem5.0}, the variety $Var\,F_0$ is unitary closed. Therefore, $F_0$ satisfies the partial linearization $f\Delta_z^m$  which has the form $\sum_J a_Jz^Jx$, where all $J$ are of the form $J=(0,j_1,\ldots,j_t)$. Without loss of generality we may assume that all $I$ in $f$ have $i_0=0$.  

Choose a lexicographically minimal $m$-tuple $I_0=(0,i_1,\ldots,i_m)$ in $f$, and 
consider a partial linearization $f\Delta$ of $f$ for
$\Delta=\Delta^{i_1}_z(z)\cdots \Delta^{i_m}_z(z^{m})$, then $f\Delta$ vanishes in $F_0$ as well (see \cite{ZSSS}). Observe
that in $F_0$ we have
$$
(z^n)R_{x,x}^m=\left\{
 \begin{array}{cc}
   0, & m>n, \\
 \tfrac{n!}{(n-m)!}\,z^{n-m}(z,x,x)^m, & m\leq n.
 \end{array}\right.
 $$
Therefore, calculating $f\Delta$ in $F_0$, we will have
$(z^Ix)\Delta=0$ for $I\neq I_0$ and $(z^{I_0}x)\Delta
=(1!)^{i_1}\cdots(m!)^{i_m}(z,x,x)^ix$. Thus
$$
0=f\Delta=\a_{I_0}(1!)^{i_1}\cdots(m!)^{i_m}(z,x,x)^ix,
$$
implying $\a_{I_0}=0$ and $f=0$.

\smallskip
Let now $f=f(z_1,\ldots,z_k,x)\in F_0[Z;x]$ be a nontrivial linear combination of basic elements of $F_0[Z;x]$ which vanishes in $F_0$ for any substitution $z_i=a_i\in K(F_0),\ x=y\in (F_0)_1$. We prove that $f=0$ by induction on $k$. Write $f$ in the form
\bes
f=\sum_{I\in \Z_+^{\infty}}f_Iz_1^I, 
\ees
where $f_I=f_I(z_2,\ldots,z_k,x)\neq 0$. Arguing as above, we may choose $I_0$ and the operator $\Delta$ of partial linearization in $z_1$ such that $(z_1^I)\Delta =0$ for $I\neq I_0$ and $(z_1^{I_0})\Delta =n(z_1,x,x)^i$ for some natural numbers $i,n$. Then 
\bes
f\Delta = f_{I_0}n(z_1,x,x)^i
\ees
is an identity in $F_0$, which implies easily that $f_{I_0}$ is an identity in $F_0$. Then by induction $f_{I_0}=0$, a contradiction. \ctd

\begin{lem}\label{lem5.2}
$Var\, (F_0)=Var\,(F_{(-1,1)}[\emptyset;x])$.
\end{lem}
\Proof
Observe first that the superalgebra $F_{(-1,1)}[\emptyset;x]$ does not satisfy identity (\ref{id4.1}). In fact, it has a base formed by the elements $(x^2)^Ix^{\varepsilon},\ I\in \Z_+^{\infty},\ \varepsilon\in\{0,1\}$. In particular, $(x^2,x,x)^k\neq 0$ in $F_{(-1,1)}[\emptyset;x]$ for any $k$. Therefore, by Proposition \ref{prop4.1}, the variety $Var\,(F_{(-1,1)}[\emptyset;x])$ is not small, that is, contains the superalgebra $F_0$. 

\smallskip

To prove the converse inclusion, it suffices to show that no  nontrivial linear combination of basic elements of the superalgebra  $F_{(-1,1)}[\emptyset;x]$ is an identity in $F_0$. Assume it is not the case, that is, $F_0$ satisfies the identity 
\bee\label{id5.1}
\sum_{I\in\Z_+^{\infty}}\a_I(x^2)^Ix^{\varepsilon}=0,\ \a_I\in\PP.
\eee
By Lemma \ref{lem5.1}, the superalgebra $F_0[z;x]$ satisfies identity 
(\ref{id5.1}) as well for any odd $x$. Let $q=\min\{d(I)\}$ for $I$ in (\ref{id5.1}),
then we have in $F_0[z;x]$ the identity $G\Delta^q_x(zx)=0$, where $G$ stands for the right part of (\ref{id5.1}) and $\Delta_x^i(y)$ is the operator of partial linearization in $x$ of degree $i$ (see \cite[ 1.4]{ZSSS}). It is clear that $((x^2)^Ix^{\varepsilon})\Delta_x^q(zx)=0$ if $d(I)>q$ (since $x^2=0$ in $F_0[z;x]$), hence we have in $F_0[z;x]$ the equality
\bee\label{id5.2}
\sum_{d(I)=q}\a_I((x^2)^Ix^{\varepsilon})\Delta^q_x(zx)=0.
\eee
Furthermore, for every $I$ with non-zero components $i_{k_1},\ldots,i_{k_q}$  we have
\bes ((x^2)^Ix^{\varepsilon})\Delta^q_x(zx)&=&((x^2)\Delta^1_x(zx))^Ix^{\varepsilon}= (x\cdot zx+zx\cdot x)^Ix^{\varepsilon}\\
&=&(\tfrac32(z,x,x))^Ix^{\varepsilon}=(\tfrac32)^q(zR_{x,x}^{k_1+1})^{i_{k_1}}\ldots,(zR_{x,x}^{k_{q}+1})^{i_{k_q}}x^{\varepsilon}.
\ees
Evidently, the obtained elements are linearly independent in $F_0[z;x]$ for different $I$, hence $\a_I=0$ for all $I$ with $d(I)=q$, and identity (\ref{id5.1}) is trivial.\ctd

\begin{lem}\label{lem5.3}
$F[Z;x]\in Var\,(F_0)$.
\end{lem}
\Proof
Again, it suffices to prove that no nontrivial linear combination $f$ of basic elements of $F[Z;x]$ is an identity in $F_0$. Write $f$ in the form
$$
f=\sum_{I\in\Z_+^{\infty}} (x^2)^If_I 
$$
where 
$$f_I=\sum_{I_1,\ldots,I_k}a_{I_1,\ldots,I_k}z_1^{I_1}\cdots z_k^{I_k}x^{\varepsilon}\neq 0,
$$
$z_i\in Z,\ I_s\in \Z_+^{\infty},\ a_{I_1,\ldots,I_k}\in\PP$. 

Let $z\in Z\setminus\{z_1,\ldots,z_k\}$ and $q=\min\{d(I)\}$ in the expression for $f$. Arguing as in the proof of Lemma \ref{lem5.2}, we have in $F_0$ the identity 
\bes
0=f\Delta^q_x(zx)=(\tfrac32)^q\sum_{I,I_1,\ldots,I_k,\, d(I)=q}a_{I_1,\ldots,I_k}(z,x,x)^Iz_1^{I_1}\cdots z_k^{I_k}x^{\varepsilon}.
\ees
By Lemma \ref{lem5.1}, this implies that all $a_{I_1,\ldots,I_k}=0$ and $f_I=0$, a contradiction. \ctd

\begin{lem}\label{lem5.4}
$Var\,F_0=Var\,F_1$.
\end{lem}
\Proof 
Since the superalgebra $F_1$ is a homomorphic image of $F_t=F[z;x]$, by Lemma \ref{lem5.3} we have $F_1\in Var\,F_0$. On the other hand,  $F_1$ evidently does not satisfy identity (\ref{id4.1}), hence  by Proposition \ref{prop4.1} it is not small and $F_0\in Var\, F_1$. \ctd  

\begin{Cor}\label{cor5.5}
$Var\, A_0= Var\,\overline{A}_0=Var\, G(F_0)$.
\end{Cor}
\Proof For a variety $\mathcal{V}$ of algebras denote by $\widetilde{\mathcal{V}}$ the variety of $\mathcal{V}$-superalgebras. Clearly, it suffices to prove that 
$\widetilde{Var\,{\overline{A}_0}}=\widetilde{Var\,A_0}=Var\,F_0$.

Recall that $\overline{A}_0$ is a homomorphic image of $A_0$, and $A_0\subseteq G(F_1)$. Therefore,
\bes
\widetilde{Var\,{\overline{A}_0}}\subseteq  \widetilde{Var\,A_0}\subseteq Var\,F_1=Var\,F_0,
\ees
and it suffices to prove that the variety $\widetilde{Var\,{\overline{A}_0}}$ is not small. Consider the superalgebra $G\otimes \overline{A}_0=G_0\otimes \overline{A}_0+G_1\otimes \overline{A}_0$ which evidently belongs to $\widetilde{Var\,{\overline{A}_0}}$. Take in $G\otimes \overline{A}_0$ the elements $z=1\otimes z$, $x_k=g_1\otimes e_1+\cdots+g_{2k}\otimes  e_{2k}$, where by $g_i$ we denote the canonical generators of the Grassmann algebra $G$, in order not confuse them with the generators $e_i$ of $\overline{A}_0$. We have
$(z,x_k,x_k)=\sum_{i\neq j}g_ig_j\otimes (z,e_i,e_j)$ and furthermore,
$(z,x_k,x_k)^k=(2k)! g_1g_2\cdots g_{2k-1}g_{2k}\otimes (z,e_1,e_2)\cdots (z,e_{2k-1},e_{2k})\neq 0$. Therefore, the superalgebra $G\otimes \overline{A}_0$ does not satisfy identity (\ref{id4.1}) for any $k$ and hence the variety $\widetilde{Var\,{\overline{A}_0}}$ is not small. \ctd

\section{Proof of Theorems A and B}

\subsection{Proof of Theorem A}  We will present the proof in  three steps.

\begin{lem}\label{lem6.1} 
$G(B) \in Var\,(\overline{A}_{0})$, where $B=B(\Gamma,D,\gamma),\ \Gamma=\Gamma_0$.
\end{lem}
\Proof The inclusion in lemma is equivalent to the inclusion $B\in \widetilde{Var\,{\overline{A}_0}}=Var\,(F_0)=Var\,(F[Z;x])$. 
Let $Z$ be such a set that $\Gamma$ is isomorphic to a homomorphic image of the polynomial ring without constant terms $\PP_0[Z]$ under a homomorphism $\varphi$.
Consider the polynomial ring $\Phi_0[T_Z]$, where 
$T_Z=T_0\cup(\cup_{z\in Z}T_z)$, $T_0=\{t_0,t_1,\ldots\},
T_z=\{z=z_0,z_1,z_2,\ldots\},\ z\in Z$, and let $D$ be the
derivation of $\Phi[T_Z]$ defined by $D(t_i)=t_{i+1},\
D(z_i)=z_{i+1}$. By Proposition \ref{prop3.6}, $F[Z;x]\cong B(\PP_0[T_Z],D,t_0)$. 
 Extend the homomorphism $\varphi$ to a homomorphism $\tilde\f:\PP_0[T_Z]\rightarrow \Gamma$ by setting 
$\tilde\f(t_i)=D^i(\gamma),\, \tilde\f(z_i)=D^i(\f(z)),\,z\in Z$; then evidently $\tilde\f$ is a homomorphism of differential algebras
$\tilde\f:(\PP[T_Z],D)\rightarrow (\Gamma,D).$ Clearly, $\tilde\f$ induces a surjective homomorphism of superalgebras
$\tilde{\tilde\f}: B(\PP_0[T_Z],D,t_0)\rightarrow B(\Gamma,D,\gamma)$. Therefore, $B(\Gamma,D,\gamma)\in Var\,(F[Z;x]).$ \ctd

\begin{lem}\label{lem6.2} 
$Var\,(G(B)) = Var\,(\overline{A}_{0})$, where $B$ is a simple non-associative (-1,1)-superalgebra.
\end{lem}
\Proof
In view of \cite{Sh2}, a simple (-1,1)-superalgebra $B=A+M,\ (A=B_0,\,M=B_1)$ over a field of
characteristic $0$ has the following properties:

\smallskip

a) $A$ is a unital differentially simple associative and commutative algebra without zero divisors and $M$ is an associative commutative $A$-bimodule;

b) If $0\neq x\in M$, then $R_{x,x}$ is a nonzero derivation of $A$.

\smallskip

It follows from a) and b) that $B$ does not satisfy identity (\ref{id4.1}) for any $k$, hence $F_0\in Var\,B$ and $ \overline{A}_{0}\subseteq G(F_0)\in Var\,(G(B))$.

\smallskip
On the other hand, the arguments from \cite{ZhSh} given for simple Jordan superalgebras with associative even part 
are applied to simple (-1.1)-superalgebras as well and show that 
the superalgebra $B=A+M$ can be embedded into a twisted superalgebra $B(\Gamma ,D,\gamma)$ of  vector type, 
where $\Gamma $ is the field of fractions of $A$. Therefore, by Lemma \ref{lem6.1}, $Var\,(G(B))\subseteq Var\,(G(B(\Gamma,D,\gamma)))\subseteq \,Var\, \overline{A}_{0}$. \ctd

\begin{lem} \label{lem6.3} Let $A_{VF}=A(V,\PP,\tau,\lambda)$ be a (-1,1)-superalgebra of vector fields over a
field $\Phi$ of characteristic $0$. Suppose that there exists $u\in V$ with $\tau(u)\neq 0$. Then  $Var\,(G(A_{VF}))=Var\,(\overline{A}_{0})$.
\end{lem}
\Proof 
By \cite[Theorem 1]{Sh1}, the superalgebra $A_{VF}$ is  of even vector type,
hence  $G(A_{VF}) \in Var\,(\overline{A}_{0})$ by Lemma \ref{lem6.1}. 
On the other hand, we have in $A_{VF}$ the equality $(a_{u},x_{v},x_{v}) =4\tau (u) a_{\lambda (u+2v)}$ which implies that 
$(a_{u},x_{v},x_{v})^k=4^k (\tau(u))^k a_{k\lambda(u+2v)}\neq 0$ for any $k$.
Therefore, by Proposition \ref{prop4.1} $F_0\in Var\,(A_{VF})$ and 
$\overline{A}_{0}\in Var\,(G(F_0)) \subseteq Var\,(G(A_{VF}))$.  \ctd

\smallskip

Theorem A now follows from Lemmas \ref{lem6.1} - \ref{lem6.3}, Lemma \ref{lem5.2}, and Corollary \ref{cor5.5}.

\smallskip

It follows from Lemma \ref{lem6.3} that the prime (-1,1)-monsters of characteristic $0$,
constructed in \cite{Pch2} and \cite{Sh1} are isomorphic.

\begin{Cor}\label{cor6.4} The free algebra of countable rank in the
variety generated by the Grassmann envelope of the (-1,1)-superalgebra of vector type $B(\Phi[t],\frac{d}{dt},0)$ over a field $\Phi$ of characteristic $0$ is prime.
\end{Cor}

\subsection{Proof of Theorem B}

The proof of Theorem B we also give in three steps.

\begin{lem} \label{lem6.5} The variety $\mathcal{V}_0:=Var\,(\overline{A}_{0})$ over a field $\Phi$ of characteristic 0 is the smallest variety of (-1,1)-algebras that contains a prime non-associative algebra.
\end{lem}
\Proof Let $\mathcal{V}$ be a variety containing a prime nonassociative algebra $A$.
Since $A$ is a strongly (-1,1)-algebra \cite{Hen1}, we may assume, without loss of generality, that $\mathcal{V}$ is strongly (-1,1). Assume that $\mathcal{V}$ does not contain  $\overline{A}_{0}$, then the corresponding variety $\tilde{\mathcal{V}}$ of $\mathcal{V}$-superalgebras does not contain $F_0$. By Proposition \ref{prop4.1},  the identity $(z,x,x)^k=0$ holds in $\tilde{\mathcal{V}}$ for any central $z$ and odd $x$.

Consider the $\mathcal{V}$-superalgebra $G\otimes A=G_0\otimes A+G_1\otimes A$. Let $z_1,\ldots,z_k\in K(A)$ and $a_1,\ldots,a_{2k}\in A$; consider
in $G\otimes A$ the elements $z=e_1e_2\otimes z_1+\cdots + e_{2k-1}e_{2k}\otimes z_k$ and $x=e_{2k+1}\otimes a_1+\cdots +e_{4k}\otimes a_{2k}$, where $e_1,\ldots,e_{4k}$ are the elements of the canonical base of $G$. 
Then $z\in K(G\otimes A),\ x\in (G\otimes A)_1$, hence we have 
\bes
0&=&(z,x,x)^k=\prod_{i,s,t} e_{2i-1}e_{2i}e_se_t\otimes (z_i,a_s,a_t)\\
&=&\pm e_1\cdots e_{4k}\otimes \prod_{i,s,t} (z_i,a_s,a_t),
\ees
where $i\in\{1,\cdots,k\},\ s,t\in\{1,\ldots,2k\}$.
In view of identities (\ref{10}), (\ref{11}), the product $(z,a,b)(z',c,d)$ for $z,z'\in K(A)$ is skewsymmetric in $a,b,c,d$. Moreover, the associators $(z_i,a_s,a_t)$ lie in the associative and commutative algebra $K(A)$. Therefore, we have 
\bes
0=(4k)!\, e_1e_2\cdots e_{4k}\otimes (z_1,a_1,a_2)\cdots (z_k,a_{2k-1},a_{2k}).
\ees
Returning to the algebra $A$, we have in it the equality $(K_{1})^{k}=0$, where $K_{1}=(K(A),A,A)$. Since the center $K(A)$ is closed under the operators $R_{x,y}$, the set $K_{1}A^{\sharp}$
is an ideal of $A$, and by induction it is easy to see that
$(K_{1}A^{\sharp})^{N}\subseteq K_{1}^{N}{A}^{\sharp}$.
Then the ideal $K_{1}{A}^{\sharp}$ is nilpotent, hence $K_{1}=0$ and $K(A) \subseteq Z(A)$. In particular, the  center $Z(A)$
contains nilpotent elements, such as commutators \cite{Kl},  a contradiction. \ctd

\begin{lem}\label{lem6.6} 
Let $F$ be the free algebra of countable rank in the variety $\mathcal{V}_0=Var\,(A_0)$. Then every  metaideal $I$ of finite index of the algebra $F$ generates the same  variety $\mathcal{V}_0$.
\end{lem}
\Proof Let $I$ be a metaideal of the algebra $F$. Since an ideal of a prime
(-1,1)-algebra is a prime algebra \cite{Roo2},  $I$ is a prime algebra. In
addition, an ideal of a non-associative prime (-1,1)-algebra can not be an
associative algebra \cite{Pch8}. Hence by Lemma \ref{lem6.5} and Corollary \ref{cor5.5}   $\mathcal{V}_0 \subseteq Var\,(I)$. The converse inclusion is evident. \ctd

\smallskip

Recall that the ideal of an algebra $A$ generated by all commutators $[a,b],\,a,b\in A$, is called the {\em commutant} of the algebra $A$.
We will denote the commutant of an algebra $A$ by $A'$.
\begin{lem}\label{lem6.7}
A metaideal of the algebra $F$ containing in its commutant $F^{\prime}$ can not be a free algebra of any variety. 
\end{lem}
\Proof 
 Let $I$ be a metaideal of $F$ containing in  $F^{\prime}$. Assume that it is a free algebra with a set of free generators $f_{1},f_{2},\ldots$. Since the commutant $F^{\prime}$ is a nil-algebra \cite{Hen1}, there exists a number $n$ such that $f_{1}^{n}=0$.
Hence $I$ is a nil-algebra of bounded index, and by \cite{Roo3} it is solvable. This is a contradiction, since a prime (-1,1)-algebra
can not be solvable. \ctd

Theorem B now follows from Lemmas \ref{lem6.5} - \ref{lem6.7}.

\section{Jordan algebras of type $A^{(+)}$ for (-1,1)-algebra $A$}

\subsection{The operator $f^{+}$} 

In a Jordan algebra $J$,  we denote by $D_{a,b}$ the inner derivation operator $D_{a,b}=\left[R_{a},R_{b}\right] :x\mapsto \left( a,x,b\right) $. If $J=A^{(+)}$, then for $a\in A$ we denote by $R_{a}^{+}$ the operator of right multiplication in the algebra $A^{(+)}$.

\smallskip
More generally, let  $f=f( x_{1},...,x_{n})$ be a nonassociative polynomial,
i.e., an element of the free nonassociative algebra $\Phi \left\{ X\right\} $. Consider the subalgebra $\Phi ^{+}\left\{ X\right\} $ of the algebra $\Phi \left\{
X\right\} ^{(+)}$ generated by the set $X$, and the homomorphism $\varphi
:\Phi \left\{ X\right\} \rightarrow \Phi ^{+}\left\{ X\right\} $ extending
the identity mapping $X$ onto itself. We set $f^{+}=\varphi (f) $. In other words, $f^{+}$ denotes the  polynomial obtained from $f$ by  replacement the multiplication $\cdot$ by its symmetrization "$\odot $". Similar notation is used for  operators. Thus, if $\rho $ is an element of the multiplication algebra of the algebra $\Phi \left\{ X\right\} $, the operator $\rho ^{+}$ is defined by $x\rho
^{+}=\left( x\rho \right) ^{+}$, where $x\in X$. Observe that this agrees with the previous notation $R_a^{+}$.

Since a right alternative algebra $A$ satisfies identity \cite[p. 69]{ZSSS}
$$
4\left( x,y,z\right) ^{+}=2\left( y,x,z\right) +\left[ y,\left[ x,z\right]
\right],
$$
every strongly (-1,1)-algebra satisfies the identity
\bee
2(x,y,z) ^{+}=(y,x,z),\label{7}
\eee
or $R_{a,b}=2D_{a,b}^{+}$, where $D_{a,b}^{+}=[R_{a}^{+},R_{b}^{+}]$ .

\subsection{The alternative center of a Jordan superalgebra} 

Following E.\,Zel\-ma\-nov \cite{Zel},  define {\em the  alternative center $Z_{alt}(J)$ of a Jordan superalgebra $J$} via
$$
Z_{alt}\left( J\right) =\{ z\in J\,|\,(z,x,y)+(-1)^{xy}(z,y,x)=0\}
$$
for any homogeneous $x,y\in J$. 
For the sake of brevity, we will call the elements of $Z_{alt}(J)$     \textit{central}. 
It is easy to see that the center $Z_{alt}\left( J\right) $ has properties similar to those of  $K\left( A\right)$. 
\begin{lem}\label{lem7.0}
The center $Z_{alt}=Z_{alt}\,(J)$ satisfies the following properties:
\begin{itemize}
\item[a)] $(Z_{alt},Z_{alt},J)=(Z_{alt},J,Z_{alt})=0$;

\item[b)] for any $z\in Z_{alt},\ a,b\in J$ holds 
$(a,z,b)=2(-1)^{az}(z,a,b)$;
\item[c)] $(Z_{alt},J,J)+(J,Z_{alt},J)\subseteq Z_{alt}$;
\item[d)]  $Z_{alt}(J)$ is a subsuperalgebra of $J$;
\item[e)] if  $A$ is a (-1,1)-superalgebra then $K(A) \subseteq Z_{alt}( A^{(+)})$.
\item[f)] if $J=J(\Gamma,D)$ is a superalgebra of vector type then $Z_{alt}(J)=\Gamma$.
\end{itemize}
\end{lem}
\Proof  Making a standart passage to Grassmann envelopes, we reduce the proof to the case of algebras.
One can directly check that any commutative algebra satisfies  the associator Jacoby identity  (\ref{6}) and also the flexibility identity
\bee
(x,y,x)&=&0\label{5}.
\eee
Let $z\in Z_{alt},\ a,b\in J$, then by  (\ref{6}) we have
$$
(z,a,b)+(a,b,z)+(b,z,a)=0.
$$
By the definition of $Z_{alt}$, $(z,a,b)=-(z,b,a)=(a,b,z)$, hence $2(z,a,b)+(b,z,a)=0$ and by linearized (\ref{5})  $2(z,a,b)=-(b,z,a)=(a,z,b)$, which proves b).
 
Now, for $z_1,z_2\in Z_{alt}$ we have by b)
$$
(z_1,z_2,a)=2(z_2,z_1,a)=4(z_1,z_2,a),
$$
hence $(z_1,z_2,a)=0$, which implies a). 

Furthermore, one can easily check that $Z_{alt}$ is invariant under derivations of $J$. In particular, $(a,Z_{alt},b)=Z_{alt}D_{a,b}\subseteq Z_{alt}$, which in view of b) implies c). Moreover, the fact that the applications $R_{a,b}$ are derivations of $Z_{alt}$ implies easily d). 

Finally, e) follows from relations (\ref{7}),(\ref{8}), and f) is proved directly.

 \ctd

\smallskip
Let us call a subspace $V\subseteq Z_{alt}(J)$ to be {\em $D$-invariant} if $(V,J,J)\subseteq V$.
\begin{Cor}\label{cor7.0}
The following properties are true.
\begin{itemize}
\item[a)]
If $U, V$ are $D$-invariant then so is $UV$;
\item[b)]
If $U$ is $D$-invariant then $UJ^{\sharp}$ is an ideal of $J$;
\item[c)]
If $U,V$ are $D$-invariant then $(UJ)(VJ)\subseteq (UV)J^{\sharp}$.  
\end{itemize} 
\end{Cor}
\Proof Let us prove c) which is the only non-evident. We have
\bes
(UJ)(VJ)&\subseteq &(UJ\cdot V)J+(UJ,V,J)=(UV\cdot J)J+(V,UJ,J)\\
&\subseteq &(UV)J+(UV,J,J)+U(V,J,J)+(V,U,J)J\\
 &\subseteq& (UV)J+UV+(UV\cdot J)J+(V\cdot UJ)J\\
&\subseteq& (UV)J^{\sharp}+(UV,J,J)+(UV\cdot J)J\subseteq   (UV)J^{\sharp}.
\ees
\ctd

\begin{prop}\label{prop7.1} In a prime nonassociative (-1,1)-algebra $A$
 the equality $K\left( A\right) =Z_{alt}\left( A^{(+)}\right)$ holds.
\end{prop}
\Proof It is known \cite{Hen1} that a prime nonassociative (-1,1)-algebra is strongly (-1,1).
If $z\in Z_{alt}\left( A^{(+)}\right) $, then, $\left( \forall x\in
A\right) \left( z,x,x\right) ^{+}=0$ and hence $z\in V\left( A\right) $ by
(\ref{7}). The same identity gives the opposite inclusion, which yields $V(A) =Z_{alt}(A^{(+)}) $. Finally, in a prime algebra $A$ we have $V\left( A\right)
=K\left( A\right) $ \cite[Lemma 21]{Pch8}. \ctd

\smallskip
 
\subsection{Functions $k(x,y;z,t)$ and $h_{x}(y,z)$.} Consider the following function in a Jordan algebra $J$:  
$$
k(x,y;z,t):=(xy,z,t) -(x,z,t)y-x(y,z,t).
$$
We will also use the notation:
$$
 k(x;y,z):=k(x,x;y,z),\ k(x;y):=k(x,x;y,y).
$$
This  function was used in the theory of Jordan algebras  by E.\,Zelmanov \cite{Zel} and V.\,Skosyrskii \cite{Skos1}. It is easy to see that if $J=A^{(+)}$ for an associative algebra $A$ then $k(x;y)=\tfrac12 [x,y]^2$ in $A$.

\begin{lem} \label{lem7.1}
The function $k(x,y;z,t)$ is symmetric in $x,y$ and in $z,t$. Moreover, $k(x,y;z,t)=k(z,t;x,y)$.
\end{lem}
\Proof It is clear that the function $k$ is symmetric in $x,y$. Furthermore,
the following identity holds in every Jordan algebra \cite{ZSSS}:
\bee\label{id7.1}
(x,yz,t)=y(x,z,t)+(x,y,z)t.
\eee
Therefore, 
\bes
k(x,y;z,t)&=&(xy,z,t)-x(y,z,t)-(x,z,t)y\\
&\stackrel{(\ref{4'})}=&-(x,y,zt)+(x,yz,t)+(x,y,z)t-y(x,z,t)\\
&\stackrel{(\ref{5}),(\ref{id7.1})}=&(zt,y,x)+y(x,z,t)+(x,y,t)z-t(z,y,x)-y(x,z,t)\\
&=&(zt,y,x)-z(t,y,x)-(z,y,x)t=k(z,t;y,x).
\ees
The obtained identitiy yields the other statements of the lemma. \ctd
\begin{lem}\label{lem7.2}
If $A$ is a strongly (-1,1)-algebra, then $k(x,y;z,t)^{+}\in K(A)$.
\end{lem}
\Proof Since the function $k(x,y;z,t)^{+}$ is symmetric in $x,y$ and in $z,t$,
it sufficies to verify that $k(x;y)^{+}\in K(A)$. Identity (\ref{4}) implies
\bee\label{id7.2}
(y,x^{2},z) =2(y,x,z) x+(y,x,[z,x]).
\eee
Therefore, in view of (\ref{7}) and (\ref{id7.2}) we have:
 \bes 
2k(x;y)^{+}&=&(y,x^2,y)-2x\odot (y,x,y)\\
&=&2(y,x,y)x+(y,x,[y,x])-2x\odot(y,x,y)\\
&=&(y,x,[y,x])+[(y,x,y),x].
\ees
Clearly, $[(y,x,y),x]\in [A,A]\subseteq K(A)$, and by (\ref{8})
$$
(y,x,[y,x])\in (A,A,K(A))\subseteq (K(A),A,A)\subseteq K(A).
$$ \ctd

\medskip

Define now the functions
\bes
Q_{a}&:&x\mapsto (a,a,x),\\
h_{x}(a,b)&:=&x[Q _{a},Q _{b}] =\left(b,b,\left( a,a,x\right)
\right) -\left( a,a,\left( b,b,x\right) \right).
\ees

\begin{lem}\label{lem7.3} If $A$ is a strongly (-1,1)-algebra, then $h_{x}^{+}(a,b) \in K(A)$.
\end{lem}
\Proof  It was proved by I.Hentzel and H.Smith \cite[identity (19)]{HenSm} that the function $h_{x}(a,b)$ is central in the variety of binary (-1,1)-algebras. If $A$ is a strongly
(-1,1)-algebra, then by (\ref{7}) we have 
$$ 
4h_{x}^{+}(a,b) =h_{x}(a,b)\in K(A).
$$ 
\ctd

\subsection{Some identities of $G(J( \Gamma ,\delta)) $}

\begin{prop} \label{prop7.4} The Grassmann envelope $G\left( J\left(
\Gamma ,\delta \right) \right) $ of a Jordan superalgebra $J(\Gamma ,\delta)$ of vector type satisfies the identities:
$$
k(x;y)Q _{c}=h_{x}(a,b)Q_{c}=k(x,y)^{2}=h_{x}(a,b)^{2}=(a,a,b)^{2}=0.
$$
\end{prop}
\Proof Consider the twisted superalgebra $B:=B(\Gamma,D,0)$ of vector type,
where $D=2\delta $. Observe that $B^{(+)}=J(\Gamma,\delta )$. The superalgebra $B$ is a strongly (-1,1), hence its Grassmann envelope $A:=G(B)$  satisfies Lemmas \ref{lem7.2} and \ref{lem7.3}, and we have the inclusions   $k(x,y)^{+},\ h_{x}^{+}\left( a,b\right) \in K\left( A\right) \subseteq Z_{alt}\left(A^{(+)}\right)$. This proves the first two equalities of the lemma. In addition, it follows from \cite[Lemma 23]{Pch8} and the Kleinfeld identity $[x,y] ^{3}=0$ which is valid in any strongly (-1,1)-algebra  that $k(x;y)^{2}=0$ in the algebra $A^{(+)}$. Furthermore, $4((a,a,b)^{+})^2=(a,a,b)^2=0$ \cite{Hen1}. Finally, by \cite[Lemma 6]{Pch8} we have  $h_{x}(a,b)^{2}=0$ in  $A^{(+)}$. Therefore, the algebra $G(J)=G( B^{(+)})=(G(B))^{(+)}$ satisfies the identities of the lemma.   \ctd

\begin{lem}\label{lem7.5} If a Jordan algebra $J$ satisfies the identities $k(x,y;z,t)=0$ and $h_{x}(y,z) =0$, then it satisfies the identity $aQ_{b}Q _{c}=0$.
\end{lem}
\Proof Let $(a,b,c)_{\sigma }$ denotes any of the six associators obtained from
the associator $(a,b,c)$ by an arbitrary permutation of variables.
The identity $k(x,y;z,t) =0$ implies that the associator $\left( a,b,c\right) $ is a derivation on all the variables. Therefore,
 \bes
4(a,(a,b,c)_{\s},b)&=&4\,a(a,b,c)_{\s}\cdot b-4\,a\cdot (a,b,c)_{\s} b\\
&=&2\,(a^{2},b,c)_{\s} b-2\,(a,b^{2},c)_{\s} a\\
&=&(a^{2},b^{2},c)_{\s} -(a^{2},b^{2},c)_{\s} =0.
 \ees
Furthermore,
\bes
(a,aQ_b,c)=-(a,(a,b,b),c)=-(a,(a,b,c),c)\Delta^1_c(b)+(a,(a,b,c),b)=0,
\ees
where $\Delta^1_c(b)$ is the operator of  partial linearization with respect to $c$ \cite[1.4]{ZSSS}.

The application $a\mapsto aQ_b$ is a derivation. Therefore, by the previous identity,
\bes
(aQ_b,a,c)&=&-(a,aQ_b,c)-(a,a,cQ_b)+(a,a,c)Q_b\\
&=&-cQ_bQ_a+cQ_aQ_b=h_c(a,b)=0.
\ees
From the Associator Jacobi identity (\ref{6}) we have
\bes
(aQ_b,c,a)=-(c,a,aQ_b)-(a,aQ_b,c)=(aQ_b,a,c)-(a,aQ_b,c)=0.
\ees
Finally,
\bes
aQ_bQ_c=((a,b,b),c,c)=-(aQ_b,a,c)\Delta^1_a(c)+(cQ_b,a,c)=0.
\ees
\ctd

\subsection{The Jordan central extension  $J[Z;x]$} Denote by $J[Z;x]$ the quotient superalgebra of the free Jordan superalgebra $Jord\,[Z;x]$ generated by a nonempty set $Z$ of even elements and an odd element $x$, by the ideal generated by the elements of the form
$$
(z,a,b) +(-1)^{ab}(z,b,a),\text{ where } z\in Z;\ a,b\in Jord[Z;x].
$$
Clearly, $J[Z;x]$ satisfies the following universal property:\\
 {\it for any  Jordan superalgebra $J$, any odd $y\in J_1$,  and for any
mapping $\varphi :Z\rightarrow Z_{alt}(J)$ there exists a unique
homomorphism $\tilde\varphi:J[Z;x]\rightarrow J$ such that
$\tilde\varphi|_Z=\varphi,\ \tilde\varphi(x)=y.$}

The  superalgebra $J[Z;x]$ plays for Jordan algebras the same role as the  superalgebra $F[Z;x]$ does in the variety $\mathcal{S}t$ of strongly 
(-1,1)-algebras.

Consider again the polynomial ring $\PP[T_Z]$ on the set of variables $T_Z=\{z_0=z,\,z_1,\ldots,\,|\,z\in Z\}$ with the derivation $D: z_i\mapsto z_{i+1},\,i=0,1,\ldots$.

\begin{prop}\label{prop7.6} 
The  superalgebra $J[Z;x]$ is isomorphic to the  subsuperalgebra $\PP_0[T_Z]\oplus\overline{\PP[T_Z]}$ of the superalgebra of vector type $J(\PP[T_Z],D)$.
\end{prop} 
 \Proof Denote $J=J[Z;x]$.
Consider  the subalgebra $A$ of $J_0$ generated by the elements $zR_{x,x}^i,\ z\in Z,\,i\geq 0$. Denote $D=R_{x,x}$, then by Lemma \ref{lem7.0} $D$ is a derivation of $Z_{alt}(J)$ and $A\subseteq Z_{alt}(J)$. Moreover, $D(A)\subseteq A$ hence $D$ is a derivation of $A$. Let us prove that $J=A+A^{\sharp}x$. Since $Z\cup\{x\}\subseteq A$, it suffices to prove that $A+A^{\sharp}x$ is a subsuperalgebra of $J$. By Lemma \ref{lem7.0} $A^{\sharp}x$ is an associative bimodule over $A$, hence we need only to consider the product $(A^{\sharp}x)(A^{\sharp}x)$. 
Let $a,b\in A^{\sharp}$, then  by Lemma \ref{lem7.0} we have
 \bes
(ax)(bx)&=&(ax\cdot b)x-(ax,b,x)=((ab)x)x-2(b,ax,x)\\
&=&(ab,x,x)-2a(b,x,x)=(ab)R_{x,x}-2a(bR_{x,x})\\
&=&D(a)b-aD(b)\in A.
\ees
Therefore, $A+A^{\sharp}x$ is a subsuperalgebra of $J$ and $J=A+A^{\sharp}x$.  Moreover, the obtained equality shows that  $J^{\sharp}$ is a homomorphic image of the superalgebra of vector type $J(A^{\sharp},D)$ under the homomorphism $a+\bar b\mapsto a+ax$. 

Consider the homomorphism $\varphi :\Phi[T_Z]\rightarrow A^{\sharp},\  z_i\mapsto D^i(z),\,z\in Z$. Clearly, it is a homomorphism of differential
algebras which can be extended to a homomorphism of superalgebras
$$
\tilde\varphi: J(\Phi[T_Z], D)\rightarrow J(A^{\sharp}, D).
$$
Now the composition $\pi\circ\tilde\varphi$ maps surjectively
$\PP_0[T_Z]\oplus\overline{\PP[T_Z]}$  onto $J$. By the universal property of $J$, in view of Lemma \ref{lem7.0}.f),  this map is invertible and hence it is an isomorphism.
 \ctd

\begin{Cor}\label{cor7.7} 
The superalgebra $J[Z;x]$ has a base of the form $B\cup Bx$ where $B$ consists of associative and commutative monomials on the elements $zR_{x,x}^i, \ z\in Z,\ i\geq 0.$
\end{Cor} 

Recall that for  $I=(i_0,i_1,\ldots,i_k)\in \Z_+^{\infty}$ we denote   
$$
a^I=a^{i_0}(aR_{x,x})^{i_1}(aR_{x,x}^2)^{i_2}\cdots (aR_{x,x}^k)^{i_k}.
$$
\begin{Cor}\label{cor7.8} 
The superalgebra $J[Z\cup\{s\};x]$ is isomorphic to the superalgebra $F[Z;x]^{(+)}$ under the isomorphism 
$$
s^{I_0}z_1^{I_1}\ldots z_k^{I_k}x^{\varepsilon} \mapsto (x^2)^{I_0}z_1^{I_1}\ldots z_k^{I_k}x^{\varepsilon}.
$$
In particular, the superalgebra $J[z;x]$ is isomorphic to the superalgebra $F[\emptyset;x]^{(+)}$.
\end{Cor} 
 \Proof By Propositions \ref{prop7.6} and \ref{prop3.6}, we have the isomorphisms
 $$
J[Z\cup\{s\};x]^{\sharp}\cong J(\PP[T_{Z\cup\{s\}}],D)\cong B(\PP[T_{Z\cup\{s\}}],D,s)^{(+)}\cong (F[Z;x]^{\sharp})^{(+)}.
$$ 
The restriction of these isomorphisms on "nonunital parts" gives the isomorphism of the corollary. \ctd

\section{Proof of Theorem C}

\subsection{Superalgebra $JF_0$ and variety $J\mathcal{V}_0$}

Following the case of $(-1,1)$-algebras, consider the superalgebra
$$
JF_0:=J[z;x]/(((z,x,x),x,x)).
$$
It is easy to see that $JF_0\cong J_0(\PP[z,s],s\tfrac{d}{dt})\cong F_0^{(+)}$ and has a base $z^n(z,x,x)^kx^{\varepsilon},\  n,k\geq 0,\ \varepsilon+n+k>0$.

\smallskip

 Similarly to Proposition \ref{prop4.1}, we have the following characterization of varieties of Jordan superalgebras that do not contain $JF_0$.

\begin{prop}\label{prop8.1} 
A variety $\mathcal{V}$ of Jordan superalgebras does not contain $JF_0$ if and only if it satisfies the identity
 \bes
(z,x,x)^k=0
 \ees
for any central $z$ and odd $x$ and some $k>0$.
\end{prop}

The proof repeats that of Proposition \ref{prop4.1} with evident modifications. 

\ctd

Denote by $J\mathcal{V}_0$ the variety of Jordan algebras generated by the Grassmann envelope $G(JF_0)$ of the superalgebra $JF_0$.

\begin{prop}\label{prop8.2}
Variety $J\mathcal{V}_0$ is generated by any one of the following algebras
\bes
&&\bar{J_0},\ A_0^{(+)},\ G_{(-1,1)}^{(+)},\\
&&  G(J[Z;x]) \hbox{ with } Z\neq\emptyset,\\
&& G(J(\Gamma,\delta)), \hbox{   where $\Gamma=\Gamma_0$ and $\Gamma^{\delta}$ is not nilpotent.}
 \ees
\end{prop}
\Proof
Observe first that  if  $Var\,A= Var\,B$  then $Var\,A^{(+)}= Var\,B^{(+)}$. Furthermore, it is easy to check that $G(A^{(+)})\cong (G(A))^{(+)}$. By the results of Sections 5 and 6, all the algebras 
\bes
\bar{A_0},\ A_0,\ G_{(-1,1)},\ G(F_{(-1,1)}[\emptyset;x]),\  G(F[Z;x]),\ G(B(\Gamma,D,\gamma))
\ees
generate the same variety $\mathcal{V}_0=Var\,A_0$. Therefore, 
the algebras $\bar{J_0}=(\bar{A_0})^{(+)}$, $A_0^{(+)},\ G_{(-1,1)}^{(+)}$,  $G(J[Z;x])\cong (G(F[Z\setminus\{s\};x]))^{(+)}$ for some element $s\in Z$, $G(J(\Gamma,\delta))\cong (G(B(\Gamma,D,0))^{(+)}$ generate the same variety $J\mathcal{V}_0$. \ctd

\smallskip
\begin{Cor}\label{cor8.3} The free algebra of countable rank in the variety $J\mathcal{V}_0$ over a field $\Phi$ of characteristic $0$ is prime and degenerate.
\end{Cor}

This was proved in \cite{Pch2}, \cite{MZ}, and in \cite{Sh1} for 
the varieties $Var\,(A_0)^{(+)}$, $Var\,J(\PP[t],\tfrac{d}{dt}),$ and $Var\,J_{VT}$, respectively, which in fact are all equal to  $J\mathcal{V}_0$.

\subsection{Minimality of variety $J\mathcal{V}_0$} 

\begin{prop}\label{prop8.4}
Let $\mathcal{V}$ be a subvariety of $J\mathcal{V}_0$ which contains a prime Jordan algebra $J$ which is not associative. Then $\mathcal{V}=J\mathcal{V}_0$.   
\end{prop}
\Proof
 Assume that $\mathcal{V}$ does not contain the algebra $\overline{J}_{0}$. Then the arguments from the proof of Lemma \ref{lem6.5} in view of Proposition \ref{prop8.1} show that the prime algebra $J$ satisfies  the equality $Z_{1}^{m}=0$ for some $m$, where $Z_{1}=(Z_{alt}(J),J,J)$.  Since the center $Z_{alt}\left( J\right) $ is invariant under the operators $R_{x,y}$, by Corollary \ref{cor7.0} the subspace $Z_{1}J^{\sharp}$  is an ideal of $J$. By induction, using Corollary \ref{cor7.0} again, it is easy to see that $(Z_{1}J^{\sharp})^{k}\subseteq Z_{1}^{k}J^{\sharp}$.
Then the ideal $Z_{1}J^{\sharp}$ is nilpotent. Let $Z_{1}^k=0$, $Z_{1}^{k-1}\neq 0$, then $(Z_{1})^{k-1}J^{\sharp}$ is a non-zero trivial ideal in $J$ which should be zero. Therefore, $Z_1=0$ and $Z_{alt}(J)$ is contained in the center $Z(J)$.

Since $J\in \mathcal{V}\subset J\mathcal{V}_0 $, it satisfies the identities of $G(J(\Gamma,D))$. By Proposition \ref{prop7.4},
the functions $k(x;y)$ and $h_x(a,b)$ produce nilpotent central elements in $J$ and therefore are zero in $J$. Then by Lemma \ref{lem7.5}  
$(a,b,b)\in Z_{alt}(J)=Z(J)$ for any $a,b\in J$. By  Proposition \ref{prop7.4} again, $(a,b,b)^2=0$, hence $(a,b,b)=0$ in $J$ and $J=Z_{alt}(J)$, a contradiction.

\ctd

\smallskip

Theorem C follows from Propositions \ref{prop8.2} and \ref{prop8.4}.

 \smallskip

Observe that, contrary to the case of variety $\VV_0$, we do not know whether $J\VV_0$ lies in every variety that contains a prime degenerate algebra. Below we show that it is true for the prime degenerate algebra related with the Jordan superalgebra of Poisson brackets \cite{MZ} and for the algebra constructed by  V.\,G.\,Skosyrskii \cite{Skos2}.

\begin{prop}\label{prop8.5}
Let $\VV_{PB}$ be the variety generated by the Grassmann envelope of the Jordan algebra of free Poisson bracket and $\VV_{Skos}$ be the variety generated by the algebra of Skosyrskii \cite{Skos2}.
Then we have
$$
J\VV_0\subseteq\VV_{Skos}\subset\VV_{PB}.
$$
\end{prop}
\Proof
It was proved in \cite{MSZ} that every Jordan superalgebra of brackets (or, in other terms, every {\it Kantor double}) can be embedded into a Jordan superalgebra of Poisson brackets. In particular, this is true for the superalgebra of even vector type  $J(\Gamma,\delta)$. Therefore, 
$$
J\VV_0=Var\,G(J(\Gamma,D))\subset \VV_{PB}.
$$
  The inclusion is strict since the Grassmann envelope of the superalgebra of Poisson brackets does not satisfies  Proposition \ref{prop7.4}.

The algebra of Skosyrskii has the form $J(G,D)_0$, where $G$ is the Grassmann algebra on infinite number of generators $e_1,\ldots,\,e_n,\ldots$ and $D$ is the derivation of $G$ defined by the condition $D(e_i)=e_{i+1}$. It follows from \cite{Sh3} that 
$J(G,D)_0\cong B(G,D,0)_0^{(+)}$, where $B(G,D,0)$ is a $(-1,1)$-superalgebra of vector type. It is clear that $B(G,D,0)_0$ is a nonassociative prime (-1,1)-algebra, hence $\VV_0\subseteq Var\, B(G,D,0)_0$. On the other hand, as it was mentioned above, $Var\, J(G,D)\subseteq \VV_{PB}$. Therefore, $J\VV_0\subseteq \VV_{Skos}\subset \VV_{PB}$. \ctd

It remains an open question whether $J\VV_0 = \VV_{Skos}$?

\section{Enveloping ${\VV_0}$-algebra for $\mathcal{F}_{J\mathcal{V}_0}[X]$}

In this section we relate the free algebra $\FF_{J\mathcal{V}_0}[X]$ and the strongly $(-1,1)$-algebra  $\FF_{\mathcal{V}_0}[X]$.

\begin{prop}\label{prop9.1} Let $\FF_{\VV_0}[X]$ be the free algebra of countable rank in the variety $\mathcal{V}_0$, and  let $J[X]$ be the subalgebra generated by set $X$ in the algebra $(\FF_{\VV_0}[X])^{+}$. Then $J[X]$ is isomorphic to the free algebra $\FF_{J\mathcal{V}_0}[X]$.
\end{prop}
\Proof Notice first that $J[X]$ is a relatively free algebra over the set $X$, that is, any relation $f(x_1,\ldots,x_n)=0$ in the algebra $J[X]$ is its identity. In fact, since this relation is an identity of $\FF_{\VV_0}[X]$, then ${f(a_1,\ldots,a_n)=0}$ for any $a_1,\ldots,a_n\in \FF_{\VV_0}[X]$, which obviously implies that $f$ is an identity in $J[X]$.

We now show that $J[X]$ and $\FF_{J\VV_0}[X]$ have the same identities. Assume that $f(x_1,\ldots,x_n)\neq 0$ in $J[X]$. Since $\VV_0=Var\,A_0$, there exist $a_1,\ldots,a_n\in A_0$ such that $f(a_1,\ldots,a_n)\neq 0$. But then $f$ is not an identity in $A_0^{(+)}$ and hence $f(x_1,\ldots,x_n)\neq 0$ in $\FF_{J\VV_0}[X]$. 

 Conversely, if $f(x_1,\ldots,x_n)\neq 0$ in the algebra $\FF_{J\VV_0}[X]$ then  there exist $a_1,\ldots,a_n\in A_0^{(+)}$ such that $f(a_1,\ldots,a_n)\neq 0$.
But  then $f(a_1,\ldots,a_n)\neq 0$ in the algebra $A_0$ as well, 
which means that $f(x_1,\ldots,x_n)\neq 0$ in $\FF_{\VV_0}[X]$ and in $J[X]$.

Thus, a relatively free algebra $J[X]$ and $\FF_{J\mathcal{V}_0}$ have the same identities, hence they are isomorphic. \ctd

\smallskip

It remains an open question on relation between the varieties $J\VV_0$ and $\VV_0$:
{\it is it true that every $J\VV_0$-algebra has an $\VV_0$-enveloping algebra?} 
One can show that the algebra $J[X]$ is not isomorphic  to the algebra $A^{(+)}$ for any strongly $(-1,1)$-algebra $A$. 

\smallskip

We recall also the following open question which was first formulated in \cite{Pch}:
{\it Is it true that $\VV_0=\mathcal{S}t$?}

\section{Acknowledgements}
The main part of the paper was done during  S.V.Pchelintsev's visit to the University of S\~ao Paulo supported by  the FAPESP (Brazil), Grant
2012/04702--7. He is grateful to Prof. I.P.Shestakov for the invitation and hospitality, to the FAPESP for the financial support, and to the University of S\~ao Paulo  for  excellent working conditions. He acknowledges also the support by Russian Foundation for Basic Research (Grant 11-01-00938-a), and by the Program "Development of the Scientific Potential of Higher Education" (Grant 2.1.1.419). 

 I.P.Shestakov was partially supported by  the FAPESP grant 10/50347-9 and CNPq grant 3305344/2009-9. He is also grateful to Max-Plank Institute f\"ur Mathematik for hospitality and excellent working conditions.

\end{document}